\documentclass[times]{nmeauth}
\pdfoutput=1
\usepackage{moreverb}

\usepackage{amssymb,amsmath,amsfonts}
\usepackage{amsthm}
\usepackage{algorithmic}
\usepackage[ruled]{algorithm2e}
\usepackage{booktabs}
\usepackage{array}
\usepackage{subfig}
\usepackage{bm}
\usepackage{epic}
\usepackage{eepic}

\newcommand\BibTeX{{\rmfamily B\kern-.05em \textsc{i\kern-.025em b}\kern-.08em
T\kern-.1667em\lower.7ex\hbox{E}\kern-.125emX}}

\def\volumeyear{2012}

\newcommand{\R}{\mathbb{R}}
\newcommand{\Eqref}[1]{Eq.~\ref{#1}}

\newcommand{\Algref}[1]{Algorithm~\ref{#1}}

\newcommand{\mini}[1]{\underset{#1}{\text{minimize}}}

\newcommand{\st}{\text{subject to }}

\newcommand{\bh}{\bm{h}}
\newcommand{\bv}{\bm{v}}
\newcommand{\bw}{\bm{w}}
\newcommand{\mA}{\bm{A}}
\newcommand{\mI}{\bm{I}}
\newcommand{\mQ}{\bm{Q}}
\newcommand{\mR}{\bm{R}}
\newcommand{\mS}{\bm{S}}
\newcommand{\mU}{\bm{U}}
\newcommand{\mW}{\bm{W}}
\newcommand{\mY}{\bm{Y}}
\newcommand{\sL}{\mathcal{L}}
\newcommand{\sY}{\mathcal{Y}}
\newcommand{\sX}{\mathcal{X}}

\newcommand{\bu}{\bm{u}}
\newcommand{\bbu}{\overline{\bm{u}}}
\newcommand{\bbp}{\overline{p}}
\newcommand{\bbT}{\overline{T}}

\begin{document}

\runningheads{P.~G.~Constantine,  E.~T.~Phipps, T.~M.~Wildey}{UQ for network multiphysics systems}

\title{Efficient uncertainty propagation for network multiphysics systems}

\author{P.~G.~Constantine\affil{1}, E.~T.~Phipps\affil{2}\corrauth, and T.~M.~Wildey\affil{2}}

\address{\affilnum{1}Stanford University, Stanford, CA 94305
\break
\affilnum{2}Sandia National Laboratories\footnotemark[2], PO Box 5800, MS-1318, Albuquerque, NM 87185}

\corraddr{Sandia National Laboratories, Optimization and Uncertainty Quantification Department, PO Box 5800, MS-1318, Albuquerque, NM 87185.}

\footnotetext[2]{Sandia National Laboratories is a multi-program laboratory managed and operated by Sandia Corporation, a wholly owned subsidiary of Lockheed Martin Corporation, for the U.S. Department of Energy's National Nuclear Security Administration under contract DE-AC04-94AL85000.}

\begin{abstract}
  We consider a multiphysics system with multiple component models coupled
  together through network coupling interfaces, i.e., a handful of scalars. If
  each component model contains uncertainties represented by a set of
  parameters, a straightfoward uncertainty quantification (UQ) study would
  collect all uncertainties into a single set and treat the multiphysics model
  as a black box. Such an approach ignores the rich structure of the
  multiphysics system, and the combined space of uncertainties can have a large
  dimension that prohibits the use of polynomial surrogate models. We propose an
  intrusive methodology that exploits the structure of the network coupled
  multiphysics system to efficiently construct a polynomial surrogate of the
  model output as a function of uncertain inputs. Using a nonlinear elimination
  strategy, we treat the solution as a composite function: the model outputs are
  functions of the coupling terms which are functions of the uncertain
  parameters. The composite structure allows us to construct and employ a reduced
  polynomial basis that depends on the coupling terms; the basis can be
  constructed with many fewer system solves than the naive approach, which results
  in substantial computational savings. We demonstrate the method on an
  idealized model of a nuclear reactor.
\end{abstract}


\maketitle

\vspace{-6pt}

\section{Introduction}
\vspace{-2pt}

In the growing field of uncertainty quantification (UQ), researchers seek to
increase the credibility of computer simulations by modeling uncertainties in
the system inputs and measuring the resulting uncertainties in simulation
outputs. A common setup in UQ includes a physical model -- typically a partial
differential equation (PDE) -- with imprecisely characterized input data, e.g.,
boundary conditions or material properties. The uncertainty in the inputs is
modeled by a set of random variables, which induces variability in the solution
of the physical model. The PDE solution is treated as a map from model inputs to
model outputs, and one can use Monte Carlo methods to approximate statistics
like moments or density functions of the output. However, if the PDE solution is
computationally expensive, then it may not be possible to compute an ensemble
with sufficiently many members for accurate statistics. In this case, a cheaper
surrogate response surface can be trained on a few input/output evaluations,
which is subsequently sampled to approximate statistics.

The procedure just described is blind to the details of the physical model; it
only needs a well-posed map from inputs to outputs. This blindness is
particularly unsatisfying when the physical model contains rich structure. For
example, many models in engineering -- e.g., fluid-structure interaction or
turbulent combustion -- include multiple physical components that are coupled
through interfaces. Such multiphysics models are almost always computationally
expensive, and thus UQ studies require surrogate response surfaces. However, the
process of collecting the uncertainties from each component model into one set
of inputs may yield a high-dimensional input space where accurate response
surfaces are infeasible. Response surfaces based on polynomial
approximation~\cite{Xiu_Karniadakis_Askey,Ghanem_Spanos_91} become exponentially
more expensive to construct as the dimension of the input space increases; this
is related to the so-called curse of dimensionality.

In this work, we propose a method for efficiently constructing accurate
polynomial surrogates of the input/ouput map that exploits the coupled structure
of a multiphysics model. Tackling all types of multiphysics models is
unreasonable, so we focus on a model problem with particular characterstics. We
will develop the method in the context of a steady, spatially discretized model
with two physical components. These components are coupled through a handful of
scalar coupling terms, which we refer to as \emph{network coupling}. Each physical component
contains its own set of parameters representing uncertain inputs, but the nature
of the coupling is assumed to be known precisely. 

We can summarize our proposed method as follows. We first transform the network
coupled multiphysics model to a low dimensional system with a nonlinear
elimination method; this approach defines the coupling terms as implicit
functions of solutions of the physical component models. We then posit a
polynomial model of the coupling terms as functions of the random inputs that
follows the formalism of polynomial chaos methods~\cite{LeMaitre2012}; the coefficients of
the polynomial approximation are defined by a Galerkin projection of the
nonlinear residual and computed with a Newton method. Up to this point, the
method can be considered standard. 

The novelty of our proposed method lies in the approach we take to reduce the
necessary number of PDE solves when computing the Newton updates for the
Galerkin coefficients. We treat the coupling terms as a set of intermediate
dependent random variables, which creates a composite structure in PDE
solutions; the PDE solutions are functions of the coupling terms, and the
coupling terms are functions of the uncertain inputs. Following an approach
similar in spirit to~\cite{Constantine2012}, we take advantage of the composite
structure to create a new polynomial basis of the intermediate random variables
on which to project the nonlinear residual and its Jacobian. We use this new
basis to produce a modified quadrature rule with many of the weights equal to
zero. Each zero weight corresponds to a PDE solve that can be ignored.

The problem setup and methodology are similar in spirit to works published
elsewhere~\cite{Arnst:2012hy,Arnst:2012cn,Arnst:2012wx}.  These papers
investigate uncertainty quantification methods for PDEs coupled at every point
in the computational domain, and numerically construct a low-dimensional
interface between the coupled physics via truncated Karhunen-Lo\`{e}ve
expansions of each problem's solution field.  Here we study network coupled
problems where the low-dimensional interface is built into the model, and thus
smoothing of the random field input data by the PDE operator is not necessary
for efficiency.  Our approach for developing the polynomial basis in the
intermediate coupling random variables and associated modified quadrature rule
is similar to the approach here~\cite{Arnst:2012cn} as well, however we propose
a more robust approach for building the basis and quadrature rule that addresses
the inherent ill-conditioning.

The remainder of the paper proceeds as follows. In Section
\ref{sec:multiphysics_uq}, we pose the simplified model of a network coupled
multiphysics system, and we briefly discuss solution strategies. We then
introduce a set of parameters into the model that represent the uncertainties,
and we employ polynomial approximation with a Galerkin projection to approximate
the coupling terms as functions of the parameters. Section
\ref{sec:transformation} describes the reduction procedure including formulating
the PDE solutions as composite functions, building a new polynomial basis, and
computing the modified quadrature rule. In Section \ref{sec:examples}, we
demonstrate the reduction procedure on a simple composite function where its accuracy can be measured, and we then
apply the full method to a simple model inspired from nuclear engineering.

\vspace{-6pt}

\section{Multiphysics systems} 
\label{sec:multiphysics_uq}

\vspace{-2pt}

We consider coupled systems of partial differential
equations (PDEs) where each component represents a physical component. 
Such systems may be coupled in various ways. We will focus on so-called network
coupling where each component 
interacts through a fixed set of scalars; this narrow focus is related to the
methods we will use to quantify uncertainty. 
An example would be semiconductor devices coupled together in an electronic
circuit.  We note that network coupling is often a
mathematical idealization of more complex interactions invoked to simplify
either computation or analysis through lower fidelity interfaces.

Consider a
general steady-state two-component coupled system that has been discretized in
space to arrive at a coupled finite dimensional problem
\begin{equation}\label{eq:general_multiphysics}
\begin{split}
f_1(u_1,v_2) &= 0,\quad f_1,u_1\in\R^{n_1}, \quad v_2=g_2(u_2)\in\R^{m_2},\\
f_2(v_1,u_2) &= 0,\quad f_2,u_2\in\R^{n_2}, \quad v_1=g_1(u_1)\in\R^{m_1}.
\end{split}
\end{equation}
Here $f_1$ and $f_2$ represent the equations defining each component
system with corresponding intrinsic variables $u_1$ and $u_2$.  The variables
$v_1$ and $v_2$ are the coupling variables between systems, and are defined by
the interface functions $g_1$ and $g_2$. The dimensions $n_1$ and $n_2$ may be
large, but the numbers of coupling terms $m_1$ and $m_2$ are small.

There are a variety of nonlinear solution schemes that can be employed to solve
\Eqref{eq:general_multiphysics}.  The simplest are relaxation approaches that
employ some solution strategy for solving each subequation in
\Eqref{eq:general_multiphysics}, such as a nonlinear
Gauss-Seidel~\cite{Porsching69}.  These algorithms are popular when coupling
disparate physical simulations together, since different solution procedures can
be employed for each sub-problem.  The disadvantage is these types of procedures
may converge very slowly or fail to converge.

Newton-type solution methods for \Eqref{eq:general_multiphysics} generally
exhibit much better convergence properties.  One approach is to form the
augmented system:
\begin{equation}
  \begin{split}
    f_1(u_1,v_2) &= 0, \\
    f_2(v_1,u_2) &= 0, \\
    	v_1 - g_1(u_1) &= 0, \\
    v_2 - g_2(u_2) &= 0 \\
 \end{split}
\end{equation}
for all of the variables $(u_1, u_2, v_1, v_2)$.  Applying the standard
Newton's method to this system then requires solving linear systems of the form
\begin{equation}\label{eq:newton}
  \begin{bmatrix}
    \frac{\partial f_1}{\partial u_1} & 0 & 0 & \frac{\partial f_1}{\partial v_2} \\
    0 & \frac{\partial f_2}{\partial u_2} & \frac{\partial f_2}{\partial v_1} & 0 \\
    -\frac{\partial g_1}{\partial u_1} & 0 & I & 0 \\
    0 & -\frac{\partial g_2}{\partial u_2} & 0 & I \\
  \end{bmatrix}
  \begin{bmatrix}
    \Delta u_1^{(k)} \\
    \Delta u_2^{(k)} \\
    \Delta v_1^{(k)} \\
    \Delta v_2^{(k)}
  \end{bmatrix} = - 
  \begin{bmatrix}
    f_1(u_1^{(k-1)}, v_2^{(k-1)}) \\
    f_2(v_1^{(k-1)}, u_2^{(k-1)}) \\
    v_1^{(k-1)} - g_1(u_1^{(k-1)}) \\
    v_2^{(k-1)} - g_2(u_2^{(k-1)})
  \end{bmatrix}
\end{equation}
for the Newton updates. However,
developing effective solution strategies for \Eqref{eq:newton} can be
challenging.  Moreover if different simulation codes are already implemented for
each physics, coupling them together in this fashion often requires significant
modifications to the source codes.

A compromise between the two that
provides Newton convergence but supports segregated solves for each component is
nonlinear elimination.  This approach works by eliminating the intrinsic
variables $u_1$ and $u_2$ from each component, relying on the implicit function
theorem.  In particular, the equation $f_1(u_1,v_2)=0$ defines $u_1$ as an
implicit function of $v_2$.  Hence the coupling equation $v_1 - g_1(u_1(v_2)) =
0$ can be written solely as a function of $v_1$ and $v_2$.  When applied to both
components we have
\begin{equation} \label{eq:nonlinear_elimination}
  \begin{split}
    h_1(v_1,v_2) \equiv v_1 - g_1(u_1(v_2)) &= 0, \;\; \mbox{ subject to
    } \;\; f_1(u_1,v_2) = 0, \\
    h_2(v_1,v_2) \equiv v_2 - g_2(u_2(v_1)) &= 0, \;\; \mbox{ subject to
    } \;\; f_2(u_2,v_1) = 0.
  \end{split}
\end{equation}

Evaluating $h_1$ and $h_2$ for a given $(v_1,v_2)$
requires full nonlinear solves $f_1=0$ and $f_2=0$ to compute the
intermediate variables $u_1$ and $u_2$.  Applying Newton's method to
\Eqref{eq:nonlinear_elimination} requires evaluation of partial derivatives such
as $\partial h_1/\partial v_1 = I$ and $\partial h_1/\partial
v_2 = -(\partial g_1/\partial u_1)(\partial u_1/\partial v_2)$ as well.  For $u_1$
that satisfies $f_1(u_1,v_2)=0$ we have by the implicit function theorem
\begin{equation} \label{eq:ne_sens}
  \frac{\partial u_1}{\partial v_2} = -\frac{\partial f_1}{\partial u_1}^{-1}\frac{\partial f_1}{\partial v_2}.
\end{equation}
If Newton's method is used to solve $f_1=0$ given $v_2$, the same Jacobian matrix
is used to compute the sensitivities $\partial u_1/\partial v_2$ by solving a
linear system with multiple right hand sides $\partial f_1/\partial v_2$, and
thus is only reasonable if the number of coupling terms is small. 
\Algref{alg:nonlinear_elimination} summarizes the nonlinear elimination
procedure. Different solvers can be used
for $f_1=0$ and $f_2=0$, and the only additional requirement is the ability to
compute the sensitivities~\eqref{eq:ne_sens}.  Moreover, the resulting network
nonlinear system is only size $m_1+m_2$, which is much smaller than the
full Newton system above.  The disadvantage is the number of nonlinear solves of
$f_1=0$ and $f_2=0$ can be quite large.
Nevertheless, this approach is attractive in many engineering problems when
strong coupling is present but different simulation codes with different
solution properties must be coupled.
\begin{algorithm}
\caption{Nonlinear elimination approach for solving~\eqref{eq:general_multiphysics}.}
\label{alg:nonlinear_elimination}
\DontPrintSemicolon
\SetAlgoNoLine
\SetKwInput{Input}{Given}
\SetKwFor{While}{while}{do}{end}
\Input{Initial guess $v_1^{(0)}$, $v_2^{(0)}$}
$l = 0$\;
\While{not converged}{
  Solve $f_1(u_1^{(l+1)},v_2^{(l)}) = 0$ for $u_1^{(l+1)}$ \;
  Solve $f_2(v_1^{(l)},u_2^{(l+1)}) = 0$ for $u_2^{(l+1)}$ \;
  Solve $\left(\frac{\partial f_1}{\partial u_1}(u_1^{(l+1)},v_2^{(l)})\right)\left(\frac{\partial u_1}{\partial v_2}\right)^{(l+1)} = -\frac{\partial f_1}{\partial v_2}(u_1^{(l+1)},v_2^{(l)})$ for $\left(\frac{\partial u_1}{\partial v_2}\right)^{(l+1)}$ \;
  Solve $\left(\frac{\partial f_2}{\partial u_2}(v_1^{(l)},u_2^{(l+1)})\right)\left(\frac{\partial u_2}{\partial v_1}\right)^{(l+1)} = -\frac{\partial f_2}{\partial v_1}(v_1^{(l)},u_2^{(l+1)})$ for $\left(\frac{\partial u_2}{\partial v_1}\right)^{(l+1)}$ \;
  Solve $\left(\begin{smallmatrix} I & -\frac{\partial g_1}{\partial u_1}(u_1^{(l+1)})\left(\frac{\partial u_1}{\partial v_2}\right)^{(l+1)} \\ -\frac{\partial g_2}{\partial u_2}(u_2^{(l+1)})\left(\frac{\partial u_2}{\partial v_1}\right)^{(l+1)} & I \end{smallmatrix} \right) \left(\begin{smallmatrix} \Delta v_1^{(l+1)} \\ \Delta v_2^{(l+1)} \end{smallmatrix} \right) = -\left(\begin{smallmatrix} v_1^{(l)}-g_1(u_1^{(l+1)}) \\ v_2^{(l)}-g_2(u_2^{(l+1)} \end{smallmatrix}\right)$ for $\left(\begin{smallmatrix} \Delta v_1^{(l+1)} \\ \Delta v_2^{(l+1)} \end{smallmatrix} \right)$ \;
  Update $v_1^{(l+1)} = v_1^{(l)} + \Delta v_1^{(l+1)}$, $v_2^{(l+1)} = v_2^{(l)} + \Delta v_2^{(l+1)}$ \;
  $l = l+1$
}
\end{algorithm}

The nonlinear
elimination procedure described in \Algref{alg:nonlinear_elimination} has some
clear connections with certain nonoverlapping domain decomposition methods.  In
fact, most nonoverlapping domain decomposition methods perform elimination of
the subdomain (component) variables and solve a Schur complement system,
sometimes referred to as a Steklov-Poincare operator \cite{QuarterVal}, for the
interface variables.  Often, a matrix-free approach is taken which avoids
assembling the Schur complement and only requires subdomain solutions for
boundary data obtained from a Newton-Krylov procedure
\cite{GlowWheel,21281,23134,yotov2001multilevel}.
The nonlinear elimination procedure described in
\Algref{alg:nonlinear_elimination} is slightly different in the sense that the
number of coupling variables is relatively small, which allows the Schur
complement to be explicitly constructed.  

\subsection{Multiphysics systems with uncertain input data}

Next we consider a multiphysics model described by
equations~\eqref{eq:nonlinear_elimination} that contains uncertainty, and we
assume that the uncertainty can be represented by a finite collection of
independent random variables. The modified version of
equations~\eqref{eq:nonlinear_elimination} become
\begin{equation}\label{eq:stoch_general_multiphysics}
\begin{aligned}
f_1(u_1(x),v_2(x),x_1) &= 0,\quad v_2(x)=g_2(u_2(x)), \\
f_2(v_1(x),u_2(x),x_2) &= 0,\quad v_1(x)=g_1(u_1(x)),
\end{aligned}
\end{equation}
where 
\begin{equation}
x_1\in\sX_1\subseteq\mathbb{R}^{s_1},\quad 
x_2\in\sX_2\subseteq\mathbb{R}^{s_2},\quad 
x=(x_1,x_2)\in\sX=\sX_1\times\sX_2,
\end{equation}
are the random variables modeling uncertainty in the system; we assume that
$\sX_1$ and $\sX_2$ are both product spaces. Let
$\rho_1=\rho_1(x_1)$ be the separable density function for $x_1$, and let
$\rho_2=\rho_2(x_2)$ be the separable density function for $x_2$ such that all
moments exist. The density function for $x$ is then $\rho=\rho_1\rho_2$. Notice
that even though $x_1$ only directly affects $u_1$ through $f_1$ (and similarly
$x_2$ affects $u_2$ through $f_2$), both the intrinsic variables $u_1$, $u_2$ and
coupling variables $v_1$, $v_2$ must be considered as functions of all of the random
variables $x$.

In effect, we have added a set of parameters to the model to represent the
uncertainties. We assume that the model~\eqref{eq:stoch_general_multiphysics} is
well-posed for all possible values of $x$. Moreover, we assume that the model
solutions $u_1$, $u_2$ and the coupling terms $v_1$, $v_2$ are smooth functions
of $x$, which will justify our use of polynomial approximation methods. 

The goal is to approximate statistics -- such as moments or density functions --
of the solutions and coupling terms given the variability induced by the
parameters $x$. We assume that computing the solution
to~\eqref{eq:stoch_general_multiphysics} given $x$ is too costly to permit an
exhaustive Monte Carlo study. Therefore, we will construct cheaper surrogate
models of the coupling terms as functions of $x$ that will permit
exhaustive sampling studies to approximate the statistics; given a value for the
coupling terms and input parameters, computing $u_1$ and $u_2$ are merely two
more PDE solves. In particular, we
will approximate $v_1(x)$, and $v_2(x)$ as a
multivariate polynomial in $x$ built as a series of orthonormal (with respect to
$\rho$) basis polynomials. Such a construction is known in the UQ literature as
a polynomial chaos expansion~\cite{Xiu_Karniadakis_Askey,Ghanem_Spanos_91,LeMaitre:2010bm}. Once the coefficients of such a series are
determined, then computing the first and second moments of the solution
(roughly, mean and variance) are simple functions of the coefficients. 

The polynomial approximation takes the form of a truncated series with $P$
terms,
\begin{equation}
\label{eq:series}
\begin{aligned}
v_1(x) &\approx \tilde{v}_1(x) \;=\; \sum_{i=0}^P \hat{v}_1^i\psi_i(x),\\
v_2(x) &\approx \tilde{v}_2(x) \;=\; \sum_{i=0}^P \hat{v}_2^i\psi_i(x).
\end{aligned}
\end{equation}
Owing to the separability of the space $\sX$, the multivariate polynomials $\psi_i$ are constructed as tensor products of univariate polynomials orthogonal with respect to the density function for each random variable, with total order at most $N$.
Thus the number of terms $P$ is related to the degree of approximation $N$ by
$P+1={N+s_1+s_2 \choose N}$. The orthogonality of the basis polynomials
$\psi_i(x)$ can be expressed as
\begin{equation}
\int_\sX \psi_i\,\psi_j\,\rho\,dx \;=\; 
\left\{
\begin{array}{cl}
1 & \mbox{ if $i=j$,}\\
0 & \mbox{ otherwise.}
\end{array}
\right.
\end{equation}
The orthogonality of the basis implies that the coefficients of the series are
the Fourier transform and can be written as
\begin{equation} \label{eq:fourier_coeff}
\hat{v}_1^i \;=\; \int_\sX v_1\,\psi_i\,\rho\,dx, \qquad 
\hat{v}_2^i \;=\; \int_\sX v_2\,\psi_i\,\rho\,dx.
\end{equation}
These integrals can be approximated with a numerical quadrature rule, where one
need only compute the coupling terms at the quadrature points in
the space $\sX$. Such a method can use existing solvers
for~\eqref{eq:stoch_general_multiphysics} given $x$, and is therefore known as
\emph{non-intrusive}.

However, the non-intrusive approach must construct the quadrature rule on the
potentially high dimensional space $\sX$, which may result in a prohibitively large
number of solves of~\eqref{eq:stoch_general_multiphysics}. Also, the
non-intrusive approach is blind to the structure of the coupled model. In what
follows, we describe an approach that takes advantage of the coupled structure
in the context of an intrusive Galerkin technique for approximating the
coefficients in~\eqref{eq:fourier_coeff}. We will exploit the low dimensionality
of the coupling terms to build a reduced polynomial basis and reduced quadrature
rule that will improve the efficiency of the Galerkin computation.

\subsection{Galerkin}\label{sec:multi_phys_galerkin}

Next we describe a standard Galerkin approach for computing the
coefficients in the series~\eqref{eq:series}. We first formulate the
model~\eqref{eq:stoch_general_multiphysics} in the form for nonlinear elimination
\begin{equation}
\label{eq:mphys_me}
\begin{array}{rcccl}
h_1(v_1,v_2,x) &=& v_1(x)- g_1(u_1(v_2(x),x_1)) &=& 0, \\
h_2(v_1,v_2,x) &=& v_2(x)- g_2(u_2(v_1(x),x_2)) &=& 0.
\end{array}
\end{equation}
The Galerkin system of equations that defines the coefficients of the polynomial
series is constructed by first substituting the polynomial approximations
$\tilde{v}_1$ and $\tilde{v}_2$ into the residuals $h_1$ and $h_2$ and then
requiring the projection of the residuals onto the polynomial basis to be
zero. More precisely, the Galerkin system can be written
\begin{equation}
\label{eq:galerkinsys}
\begin{aligned}
\hat{h}_1^i &= \int_\sX h_1(\tilde{v}_1,\tilde{v}_2)\,\psi_i\,\rho\,dx \;=\; 0,\\
\hat{h}_2^i &= \int_\sX h_2(\tilde{v}_1,\tilde{v}_2)\,\psi_i\,\rho\,dx \;=\; 0,
\end{aligned}
\end{equation}
for $i=0,\dots,P$. Define the vectors
\begin{align}
\hat{\bv}_1 &= [\hat{v}_1^0,\dots,\hat{v}_1^P]^T,\\
\hat{\bv}_2 &= [\hat{v}_2^0,\dots,\hat{v}_2^P]^T,\\
\hat{\bh}_1 &= [\hat{h}_1^0,\dots,\hat{h}_1^P]^T,\\
\hat{\bh}_2 &= [\hat{h}_2^0,\dots,\hat{h}_2^P]^T.
\end{align}
Then the Galerkin system can be written compactly as
\begin{equation}
\begin{bmatrix}
\hat{\bh}_1\\
\hat{\bh}_2
\end{bmatrix}
\;=\;
\begin{bmatrix}
\hat{\bh}_1(\hat{\bv}_1,\hat{\bv}_2)\\
\hat{\bh}_2(\hat{\bv}_1,\hat{\bv}_2)
\end{bmatrix}
\;=\;
\begin{bmatrix}
0\\
0
\end{bmatrix}.
\end{equation}
The Newton update to the coefficients $[\delta\hat{\bv}_1,\delta\hat{\bv}_2]^T$
is computed by solving the system
\begin{equation}
\label{eq:newtongalerkin}
\begin{bmatrix}
\frac{\partial \hat{\bh}_1}{\partial \hat{\bv}_1} 
& \frac{\partial \hat{\bh}_1}{\partial \hat{\bv}_2}\\
\frac{\partial \hat{\bh}_2}{\partial \hat{\bv}_1} 
& \frac{\partial \hat{\bh}_2}{\partial \hat{\bv}_2} 
\end{bmatrix}
\begin{bmatrix}
\delta \hat{\bv}_1 \\ \delta \hat{\bv}_2
\end{bmatrix}
= 
\begin{bmatrix}
-\hat{\bh}_1 \\ -\hat{\bh}_2
\end{bmatrix}.
\end{equation}
The top left block of the Jacobian $\frac{\partial \hat{\bh}_1}{\partial
 \hat{\bv}_1}$ is equal to an identity matrix of size $(P+1)m_1$. Similarly, the
bottom right block of the Jacobian is an identity matrix of size $(P+1)m_2$. The top
right block $\frac{\partial \hat{\bh}_1}{\partial \hat{\bv}_2}$ has size
$(P+1)m_1\times (P+1)m_2$, and has a block structure within itself. For $i,j=0,\dots,P$,
the $(i,j)$ sub-block of size $m_1\times m_2$ is
\begin{align}
\frac{\partial \hat{h}_1^i}{\partial \hat{v}_2^j}
&= -\int_\sX \frac{\partial g_1}{\partial v_2}(u_1)\,\psi_i\,\psi_j\,\rho\,dx\\
&\approx -\sum_{k=0}^P \hat{g}_1^k\,\int_\sX \psi_i\,\psi_j\,\psi_k\,\rho\,dx,
\end{align}
where (with a slight abuse of notation)
\begin{equation}
\hat{g}_1^k \;=\; \int_\sX \frac{\partial g_1}{\partial v_2}(u_1)\,\psi_k\,\rho\,dx
\end{equation}
is the $k$th coefficient of a polynomial approximation of the $\partial
g_1/\partial v_2$ evaluated at $u_1$. The bottom left block of the Jacobian is
similar. 

\begin{figure}[ht]
\centering
\subfloat{
\includegraphics[width=0.45\textwidth]{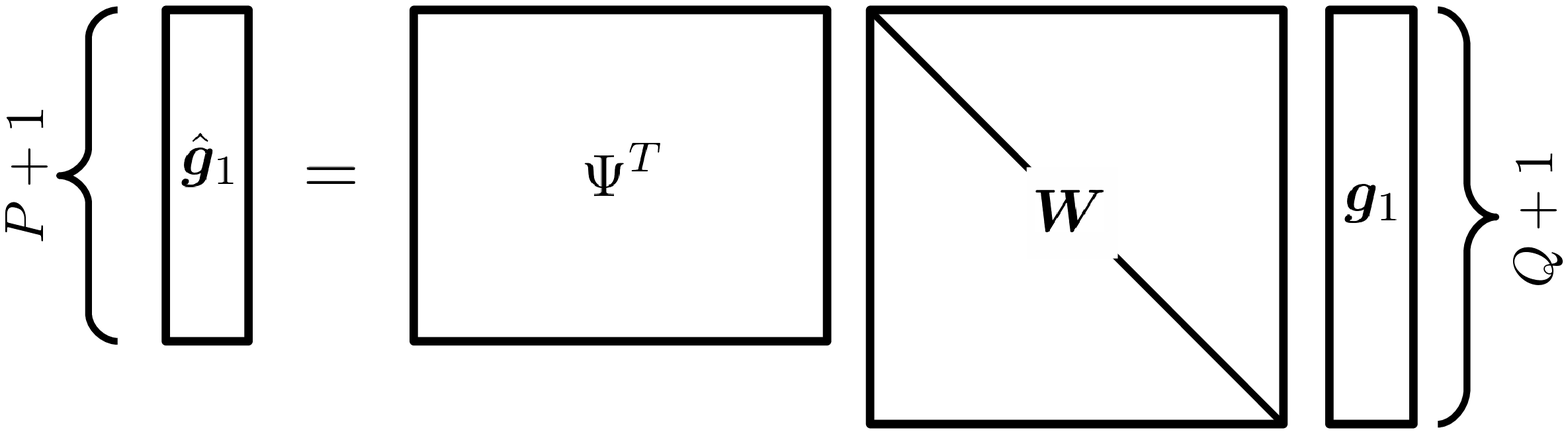}
        }\quad
\subfloat{
\includegraphics[width=0.45\textwidth]{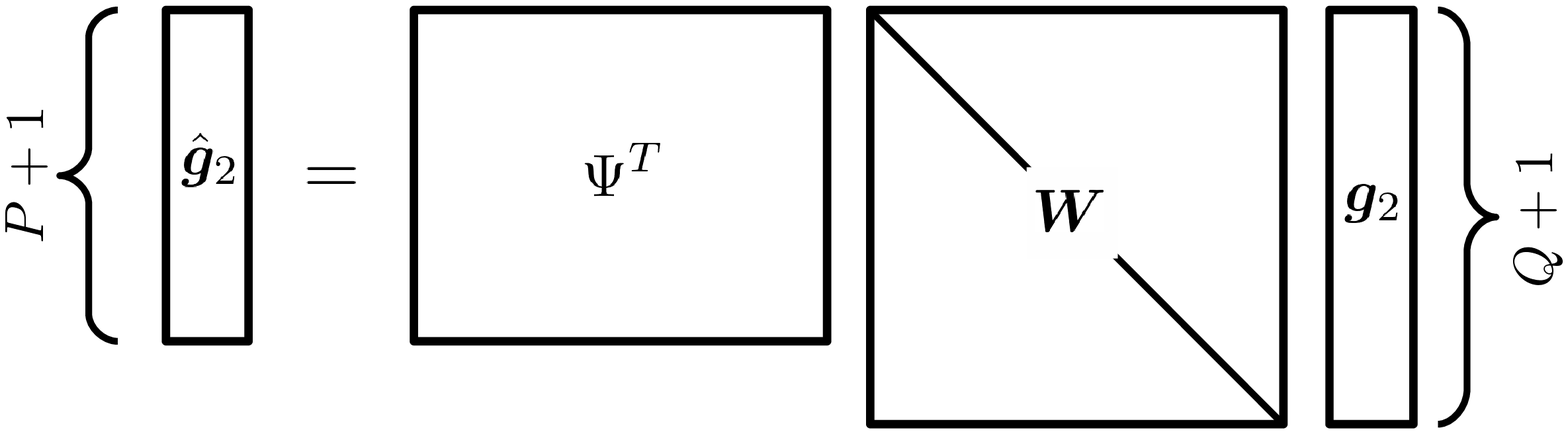}
        }\\
\subfloat{
\includegraphics[width=0.45\textwidth]{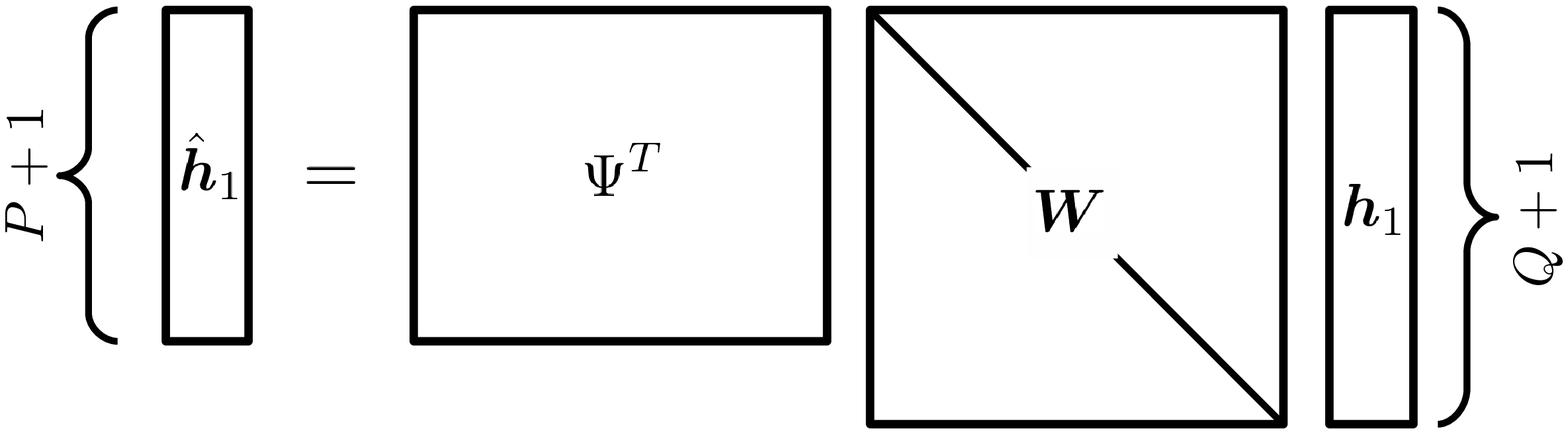}
        }\quad
\subfloat{
\includegraphics[width=0.45\textwidth]{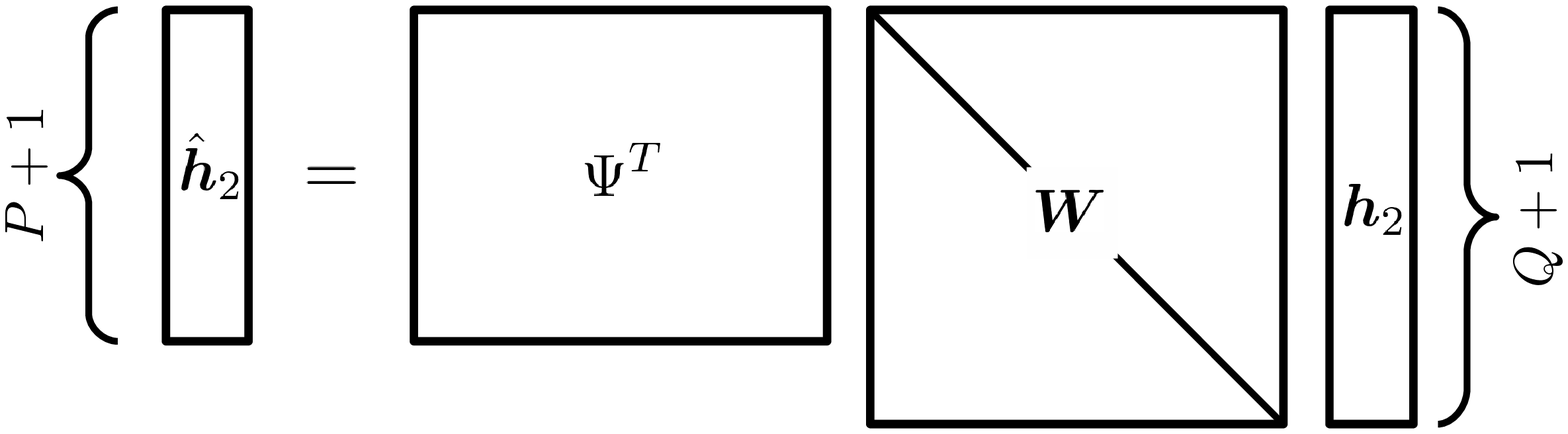}
        }
\caption{Matrix representation of pseudospectral coefficient approximation for
  the vectors of coefficients $\hat{\bm{g}}_1$, $\hat{\bm{g}}_2$,
  $\hat{\bm{h}}_1$, and $\hat{\bm{h}}_2$.}
\label{fig:ps}
\end{figure}

To construct the matrix and right hand side for the Newton update
in~\eqref{eq:newtongalerkin} we need the coefficients $\hat{g}_1^i$,
$\hat{g}_2^i$, $\hat{h}_1^i$, and $\hat{h}_2^i$. We compute these with a
pseudospectral approach, i.e., approximating the integration with a quadrature
rule. Let $\{(x^k,w^k)\}$ with $k=0,\dots,Q$ be the points and weights of a
numerical integration rule for functions defined on $\sX$; we assume the weights
are positive and associated with the measure $\rho$ such as in Gaussian
quadrature rules. Then we approximate
\begin{equation} \label{eq:ps}
\begin{aligned}
\hat{g}_1^i &\approx \sum_{k=0}^Q  \frac{\partial g_1}{\partial
  v_2}(u_1(x^k))\,\psi_i(x^k)\,w^k,\\
\hat{g}_2^i &\approx \sum_{k=0}^Q  \frac{\partial g_2}{\partial
  v_1}(u_2(x^k))\,\psi_i(x^k)\,w^k,\\ 
\hat{h}_1^i &\approx \sum_{k=0}^Q
  h_1(\tilde{v}_1(x^k),\tilde{v}_2(x^k))\,\psi_i(x^k)\,w^k,\\ 
\hat{h}_2^i &\approx \sum_{k=0}^Q
  h_2(\tilde{v}_1(x^k),\tilde{v}_2(x^k))\,\psi_i(x^k)\,w^k.
\end{aligned}
\end{equation}
These integrations require the quantities $\partial g_1/\partial v_2$, $\partial
g_2/\partial v_1$, $h_1$, and $h_2$ evaluated at the quadrature points $x^k$,
which are computed from $Q+1$ solves of the physical systems $f_1$ and $f_2$;
see~\eqref{eq:stoch_general_multiphysics}. Define the vectors
\begin{equation}
\begin{aligned}
\bm{g}_1 &= \left[
\frac{\partial g_1}{\partial v_2}(u_1(x^0)), \dots, \frac{\partial g_1}{\partial
  v_2}(u_1(x^Q))
\right]^T, \\
\bm{g}_2 &= \left[
\frac{\partial g_2}{\partial v_1}(u_2(x^0)), \dots, \frac{\partial g_2}{\partial
  v_1}(u_2(x^Q))
\right]^T, \\
\bh_1 &=
\left[h_1\big(\tilde{v}_1(x^0),\tilde{v}_2(x^0)\big),\dots,h_1\big(\tilde{v}_1(x^Q),\tilde{v}_2(x^Q)\big)\right]^T,\\ 
\bh_2 &=
\left[h_2\big(\tilde{v}_1(x^0),\tilde{v}_2(x^0)\big),\dots,h_2\big(\tilde{v}_1(x^Q),\tilde{v}_2(x^Q)\big)\right]^T.
\end{aligned}
\end{equation}
Then the pseudospectral approximations in~\eqref{eq:ps} can be written
conveniently in matrix form as in Figure \ref{fig:ps},
where $\bm{W}$ is a diagonal matrix of the quadrature weights, and $\Psi$
is a $(Q+1)\times (P+1)$ matrix of the polynomial basis evaluated at the quadrature
points, i.e., $\Psi(j,i) = \psi_i(x^j)$. 

Thus, in this standard formulation of the Galerkin method, we need $Q+1$
evaluations of the multiphysics residual and sensitivities (i.e., partial
derivatives) to construct the system for each Newton solve; this requires $Q+1$
solves of the physical components $f_1$ and $f_2$. In what follows, we
describe a procedure for constructing a reduced basis and a modified quadrature
rule that will decrease the number of residual and derivative evaluations -- and
consequently PDE solves -- needed to approximate the coefficients
in~\eqref{eq:ps}. The modified quadrature rule will have the same points in the
space $\sX$, but many of the weights will be zero; each quadrature point with a
zero weight can be ignored.

\vspace{-6pt}

\section{Reduced basis and reduced quadrature for composite functions}
\label{sec:transformation}

\vspace{-2pt}

We will develop the procedure for determining and applying the reduced basis and
modified quadrature in the context of general composite functions. This
simplifies the notation and provides a more general framework. In Section
\ref{sec:redgalerkin}, we apply this procedure to reduce the work needed to
construct the systems in the Newton solves~\eqref{eq:newtongalerkin}. 

Consider a general composite function $h(x)=h(y(x))$ for
$x\in\sX\subset\mathbb{R}^s$, where
\begin{equation}
y\;:\; \sX \rightarrow \sY\subset\mathbb{R}^m,\qquad h\;:\;\sY\rightarrow\mathbb{R}.
\end{equation}
We wish to compute the pseudospectral coefficients (i.e., discrete Fourier transform) of
$h(x)$ for the polynomial basis $\psi_i(x)$,
\begin{equation}
\begin{aligned}
\hat{h}_i &= \int_\sX h\,\psi_i\,dx \\
&\approx \sum_{j=0}^Q h(y(x^j))\, \psi_i(x^j)\,w^j.
\end{aligned}
\end{equation}
In matrix form, 
\begin{equation}
\text{\includegraphics[scale=0.4]{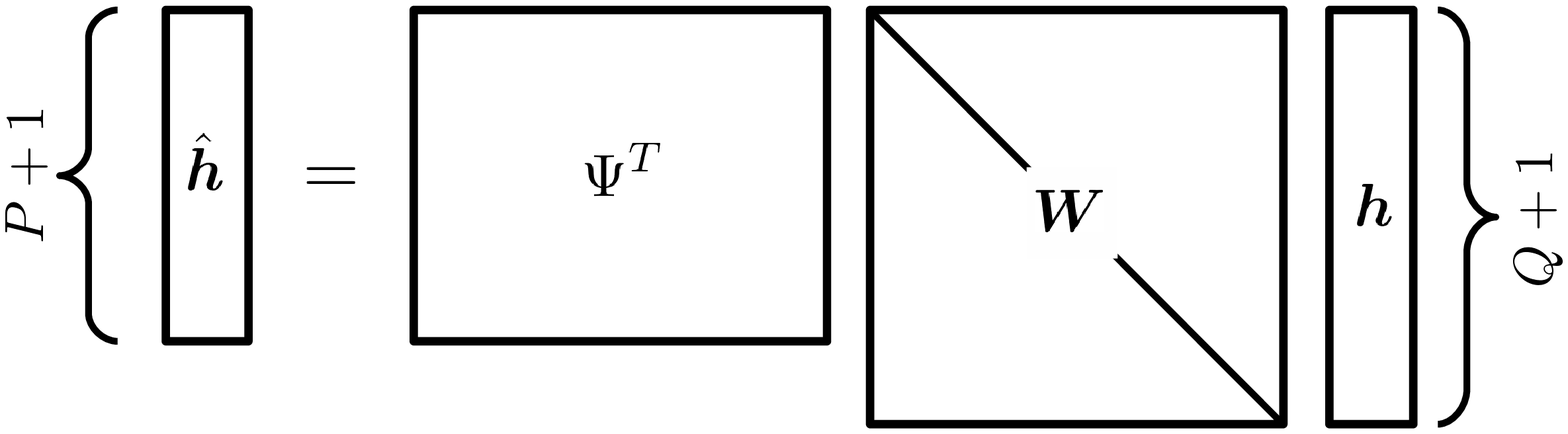}}
\end{equation}
where
\begin{equation}
\hat{\bh} = \begin{bmatrix} \hat{h}_0 \\ \vdots \\ \hat{h}_P \end{bmatrix}, \quad
\bh = \begin{bmatrix} h(y(x^0)) \\ \vdots \\ h(y(x^Q)) \end{bmatrix}, \quad
\mW=\begin{bmatrix}  w^0 & & \\ & \ddots & \\ & & w^Q \end{bmatrix},
\end{equation}
and $\Psi(j,i)=\psi_i(x^j)$ for $j=0,\dots,Q$ and $i=0,\dots,P$. This
computation requires $Q+1$ evaluations of $y(x)$ and $Q+1$ evaluations of
$h(y)$. We wish to take advantage of the composite structure to
reduce the number of evaluations of $h(y)$ in this computation.

\subsection{Reduced basis}

The first step is to construct an appropriate change of basis to approximate the
function $h(y)$. In particular, we will build a set of multivariate polynomials of the
intermediate variables $y$. We first show how to approximate $h(y)$ assuming we
have the new basis. We will then discuss how to construct the basis. 

\subsubsection{Approximation with the reduced basis}

Let $\{\phi_k(y)\}$ be a multivariate orthonormal polynomial basis in $y$ up to
degree $N'$. The number of basis elements is then $P'+1={N'+m \choose m}$ so
that $k=0,\dots,P'$; we generally assume that $P'<P$. We approximate $h$ as a
polynomial in $y$,
\begin{equation}\label{eq:pc_y}
h(y)\;\approx\; \sum_{k=0}^{P'} \hat{h}_k^y\,\phi_k(y),
\end{equation}
where the coefficients are computed as integrations with respect to the measure
of $y$, which we denote by $dy$,
\begin{equation}
\hat{h}_k^y = \int_\sY h(y) \,\phi_k(y)\,dy.
\end{equation}
The integrations can be transformed to integrals over $\sX$ and subsequently
approximated by quadrature,
\begin{equation}
\label{eq:ps_y}
\begin{aligned}
\hat{h}_k^y &= \int_\sX h(y(x)) \,\phi_k(y(x))\,dx \\
&\approx \sum_{j=0}^Q h(y(x^j))\,\phi_k(y(x^j))\,w^j.
\end{aligned}
\end{equation}
Given this representation, we can approximate the desired pseudospectral
coefficients as
\begin{equation}
\label{eq:ps_approx}
\hat{h}_i \;\approx\; \sum_{j=0}^Q \left(\sum_{k=0}^{P'} \hat{h}_k^y\,\phi_k(y(x^j))\right)\, \psi_i(x^j)\,w^j,
\end{equation}
Written compactly in matrix form, 
\begin{equation}
\label{eq:ps_approx_mat}
\text{\includegraphics[scale=0.3]{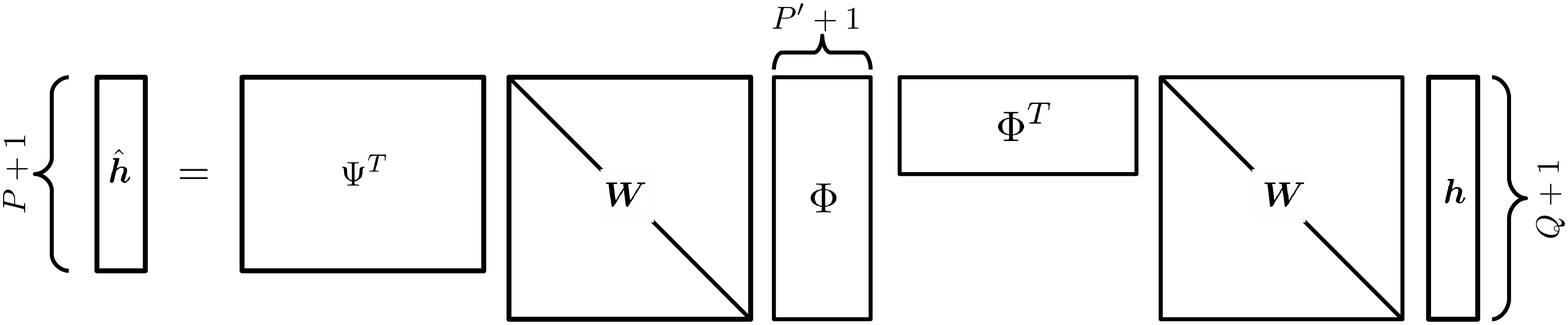}}
\end{equation}
where $\Phi(j,k)=\phi_k(y(x^j))$ for $j=0,\dots,Q$ and $k=0,\dots,P'$.

We have not yet reduced the amount of work for the approximation. Next we
describe how to construct the basis $\{\phi_k(y)\}$, and then we discuss the
modification to the quadrature rule that results in reduced computational
effort. 

\subsubsection{Constructing the reduced basis}
\label{sec:redbasis}

It is important to note that the $m$ components of $y$ are not
independent. Therefore, the multivariate polynomial bases $\phi_k(y)$ cannot be
constructed as products of univariate polynomial bases of each component of
$y$. Instead, we must use a more general procedure based on multivariate
monomials of the components of $y$. From the monomial terms, we can construct a
multivariate orthonormal basis by appropriate linear combinations, where the
coefficients are computed with a modified Gram-Schmidt approach~\cite{Bjorck1996}. 

However, in the computations~\eqref{eq:ps_y} we only need the basis polynomials
$\phi_k(y)$ evaluated at the points $y(x^j)$ corresponding to the original
numerical integration rule on $\sX$. Since we are only interested in discrete
quantities, the basis can be constructed using standard tools of numerical
linear algebra. This will also enable us to exploit standard numerical
strategies for handling the potentially unwieldy monomial terms.

We first evaluate all $m$ components of $y(x^j)$ for each quadrature point
$x^j$. We then construct a matrix of multivariate monomials from these
components. Define the set of $m$-term multi-indices $\sL$ by
\begin{equation}
\sL \;=\; \left\{(\ell_1,\dots,\ell_m)\in\mathbb{N}^m \;:\; \ell_1 + \cdots +
  \ell_m \leq N'\right\}.
\end{equation}
The number of multi-indices in $\sL$ is exactly $P'+1$, which is the number of
basis polynomials $\phi_k(y)$. The columns of the matrix of monomials will be
uniquely indexed by the multi-indices in $\sL$. For
$(\ell_1,\dots,\ell_m)\in\sL$, define the $(Q+1)\times (P'+1)$ matrix $\mY$ by 
\begin{equation}
\label{eq:monomials}
\mY(j,(\ell_1,\dots,\ell_m)) \;=\; y_1^{\ell_1}(x^j)\cdots y_m^{\ell_m}(x^j),
\end{equation}
where the superscripts $\ell$ naturally represent powers of the components of
$y$ while the superscripts $j$ represent an index to a quadrature point.
We can compute the basis $\phi_k$ evaluated at the points $y(x^j)$ with a
weighted, modified Gram-Schmidt procedure~\cite{Bjorck1996}, which yields
\begin{equation}\label{eq:basis_qr}
\mY\;=\;\Phi\mR,\quad\mbox{ such that }\quad \Phi^T\mW\Phi \;=\;\mI.
\end{equation}
The matrix $\mR$ has size $(P'+1)\times (P'+1)$ and is upper triangular. The matrix
$\Phi$ has size $(Q+1)\times (P'+1)$ and its elements are $\Phi(j,k)=\phi_k(y(x^j))$, which
exactly is what is required in~\eqref{eq:ps_approx_mat}. In practice, the
computation of $\Phi$ may benefit from one or two repetitions of the
Gram-Schmidt procedure if the matrix $\mY$ is poorly conditioned. Also, scaling
the columns of $\mY$ so that they have unit norm under the weights $w^k$ may
increase accuracy. These heuristics help ensure that the orthogonality condition
$\Phi^T\mW\Phi \;=\;\mI$ is satisfied to numerical precision. 

\subsection{Modified quadrature}
\label{sec:redquad}

From~\eqref{eq:ps_y}, it is clear that if $w^j$ is zero, then the corresponding
$h(y(x^j))$ need not be computed. This is the key motivation to seek a modified
quadrature rule. In particular, we seek a quadrature rule with the same points
$x^j$ but with modified weights such that as many of the weights as possible are
zero while retaining discrete orthogonality\footnote{Discrete orthogonality of the basis, i.e., orthogonality of the basis under the discrete inner product defined by the quadrature rule, 
  is essential to justify the pseudospectral method of computing the
  coefficients of the polynomials series.} of the reduced basis $\Phi$. We can state
such a problem formally as
\begin{equation}\label{eq:lin_prog}
\begin{array}{cl}
\mini{\bu} & \|\bu\|_0 \\
\st & \Phi^T\mU\Phi=\mI, \quad \mU=\text{diag}\,(\bu),\\
 & \bu\geq 0,
\end{array}
\end{equation}
where the zero-norm $\|\bu\|_0$ counts the number of non-zeros in $\bu$. The
problem can be reformulated as a standard linear program as follows. 
Notice that the orthogonality constraint can be rewritten as
\begin{equation}
\sum_{j=0}^Q \phi_{k_1}(y(x^j))\,\phi_{k_2}(y(x^j))\,u^j \;=\; 
\left\{
\begin{array}{lc} 1 & \mbox{ if $k_1=k_2$, } \\ 0 & \mbox{ otherwise.}\end{array}
\right.
\end{equation}
Define the $(Q+1)\times (P'+1)^2$ matrix $\mA$ by 
\begin{equation}
\mA(j,(k_1,k_2)) = \phi_{k_1}(y(x^j))\,\phi_{k_2}(y(x^j)),
\end{equation}
where the columns are uniquely indexed by the multi-indices $(k_1,k_2)$. 
Since the vector of weights $\bw=[w^0,\dots,w^Q]^T$ satisfies the orthogonality
condition by construction, we can express orthogonality condition by the linear
constraints
\begin{equation}\label{eq:red_constraint}
\mA^T\bu = \mA^T\bw.
\end{equation}
However, the matrix $\mA$ does not have full column rank, which can cause
problems in practice for linear program solvers. Therefore we seek a matrix with
the same column space of $\mA$ that is full rank; this matrix will serve as the
set of linear equality constraints for the linear program. In theory, the rank
of $\mA$ should be exactly ${2N'+m\choose m}$, since the elements
$\phi_{k_1}(y(x^j))\phi_{k_2}(y(x^j))$ are polynomials in $y$ of degree at most
$2N'$. However, in practice, we find that finite precision computations can
create a matrix $\mA$ whose rank is not exactly ${2N'+m\choose m}$. Therefore,
to find a full rank constraint set, we resort to numerical heuristics. In
particular, we use a column-pivoted, weighted, modified Gram-Schmidt procedure
as a heuristic method for a rank-revealing QR factorization~\cite{Golub1996}.  To estimate
the rank, we examine the diagonal elements of the upper triangular factor. We
compute
\begin{equation}\label{eq:red_quad_qr}
\mA\Pi = \mQ\mS, \qquad \mQ^T\mW\mQ\;=\;\mI,
\end{equation}
where $\Pi$ is a permutation matrix and $\mS$ is upper triangular.
We then examine the diagonal elements of $\mS$ and find the largest $R$ such
that $|\mS(R,R)|>\text{TOL}$ for a chosen tolerance TOL. Let $\mQ_R$ be the
first $R$ columns of $\mQ$. We assume that the parameters of the procedure
(e.g., the degree of polynomial approximations and the number of quadrature
points) is such that $R<Q+1$. Ultimately we express the linear program as
\begin{equation}\label{eq:lin_prog_final}
\begin{array}{cl}
\mini{\bu} & 0^T\bu \\
\st & \mQ_R^T\bu=\mQ_R^T\bw,\\
 & \bu\geq 0.
\end{array}
\end{equation}
The zero vector in the objective is merely to set the problem in the standard
form for a linear program solver. In fact, we only need a vector $\bu$ with
positive elements that
satisfies the linear equality constraints. 

By using a simplex method~\cite{Ye2008} to solve the linear program, we recover a vector
$\bu_\star$ with exactly $R$ non-zeros. The orthogonality with the modified
quadrature weights can be checked. If the basis $\Phi$ is not sufficiently
orthogonal with respect to the modified weights, $\bu_\star$, then the procedure
can be repeated with a stricter TOL for determining the numerical rank of $\mA$.

We have found in practice that the quality of the linear program solver can make
a significant difference in the computation of the reduced quadrature
weights. We strongly recommend solvers such as Gurobi~\cite{gurobi} that use
extended precision when dealing with the constraints; in the examples below, we
use Clp~\cite{ClpURL}.

We incorporate the modified weights into the computation in~\eqref{eq:ps_y} as
\begin{equation}
\hat{h}_k^y \;\approx\; \sum_{j=0}^Q h(y(x^j))\,\phi_k(y(x^j))\,u_\star^j.
\end{equation}
If any $u_\star^j$ is zero, then we can avoid the corresponding computation of
$h(y(x^j))$. We write the approximation in matrix form as 
\begin{equation}
\label{eq:ps_redapprox_mat}
\text{\includegraphics[scale=0.3]{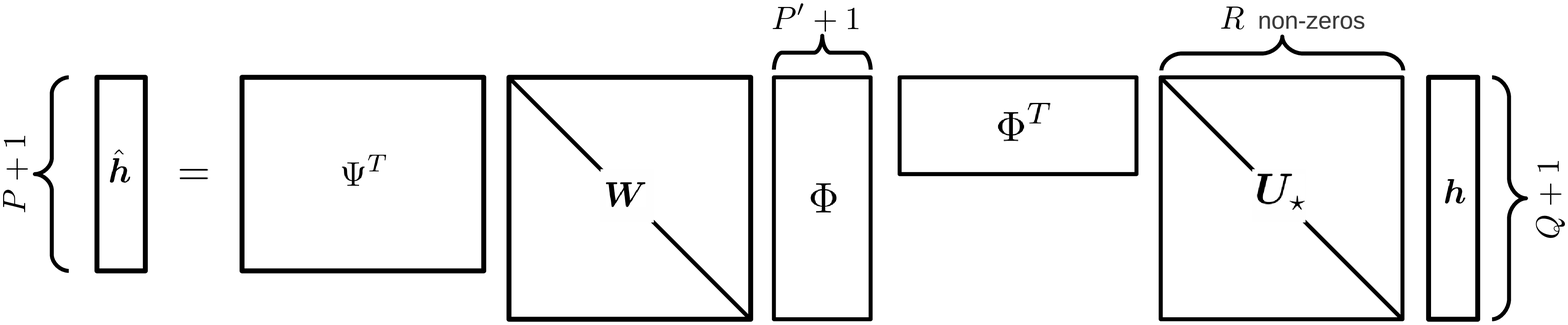}}
\end{equation}
where $\bm{U}_\star$ is a diagonal matrix with the modified weights
$\bm{u}_\star$. Each zero entry of $\bm{u}_\star$ corresponds to an evaluation
of $h(y(x))$ that can be avoided, thus saving computational effort.

\subsection{Approximating the Galerkin coefficients}
\label{sec:redgalerkin}

We next apply the reduction procedure for composite functions to the coefficients
needed to construct the Newton system~\eqref{eq:newtongalerkin}, i.e., the
computations represented in Figure \ref{fig:ps}. 

Recall that $\tilde{v}_1$ and $\tilde{v}_2$ are the series
approximations~\eqref{eq:series} to the coupling terms $v_1$ and $v_2$ in the
multiphysics model~\eqref{eq:stoch_general_multiphysics}. Consider the following
change of variables:
\begin{equation}
y_1(x) = (x_1(x),\tilde{v}_2(x)),\qquad y_2(x) = (x_2(x),\tilde{v}_1(x)).
\end{equation}
Formally, $x_1(x)=x_1$ and $x_2(x)=x_2$ select the proper subset of inputs from
the full set of original inputs; writing it this way makes the composite
structure of the problem clearer. Under this change of variables, we treat the
coupling terms $\tilde{v}_1$ and $\tilde{v}_2$ as coordinate variables. These
variables are not independent since $\tilde{v}_1=\tilde{v}_1(x_1,x_2)$ and
$\tilde{v}_2=\tilde{v}_2(x_1,x_2)$. But independence was not a prerequisite for
the reduction procedure; we avoided the need for this assumption by using the
monomials in~\eqref{eq:monomials}.

Under this change of variables, we can view the solutions $u_1$ and $u_2$ of the
coupled physical components implied in the Galerkin
system~\eqref{eq:galerkinsys} as composite functions:
\begin{equation}
  u_1=u_1(x_1,\tilde{v}_2)=u_1(y_1(x)),\qquad u_2=u_2(x_2,\tilde{v}_1)=u_2(y_2(x)).
\end{equation}
Given $y_1$ and $y_2$, computing $u_1$ and $u_2$ requires solving an expensive
PDE. However, evaluating $y_1$ and $y_2$ as functions of $x$ is relatively
cheap. To see this, first note that given $x$, evaluating $x_1$ and $x_2$ is
trivial. Second, recall that at any point in the Newton procedure we have
estimates of the coefficients $\hat{v}_1^i$ and $\hat{v}_2^i$. With these
coefficients, we can compute $\tilde{v}_1$ and $\tilde{v}_2$ using the
series~\eqref{eq:series}. Thus, we are in the situation where the reduction
procedure for composite functions will be most beneficial. Specifically, we will
reduce the number of PDE solves needed to approximate the coefficients
represented by Figure \ref{fig:ps}.

For each of the intermediate variables $y_1$ and $y_2$, we construct separate
polynomial basis sets evaluated at the quadrature points $x^j$ as described in
Section \ref{sec:redbasis}. Denote these $(Q+1)\times (P'+1)$ matrices by $\Phi_1$ and
$\Phi_2$. For each basis set, we follow the procedure in Section
\ref{sec:redquad} for computing a modified set of quadrature weights with $R_1$ and $R_2$
non-zeros respectively. Denote these weight sets by the vectors $\bu_\star^1$ and
$\bu_\star^2$. 

\begin{figure}[ht]
\centering
\subfloat{
\includegraphics[width=0.9\textwidth]{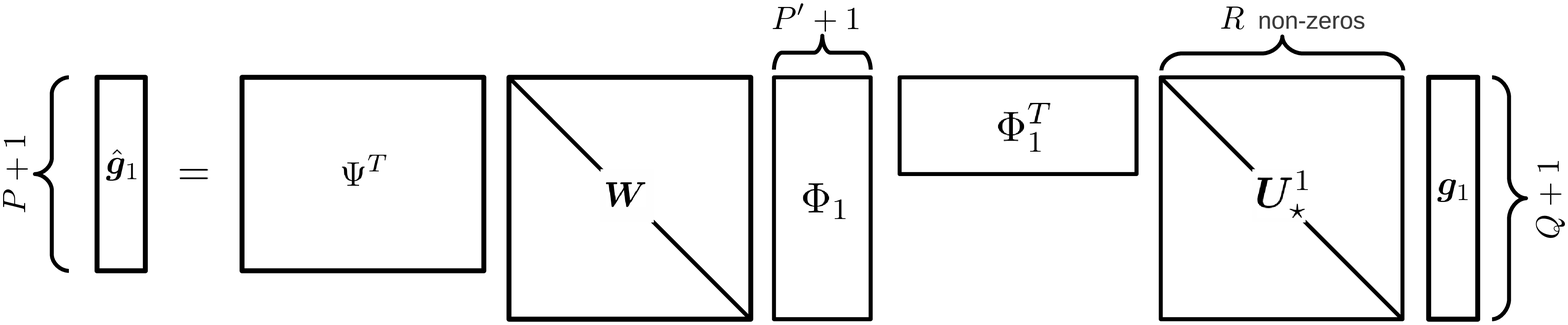}
        }\\
\subfloat{
\includegraphics[width=0.9\textwidth]{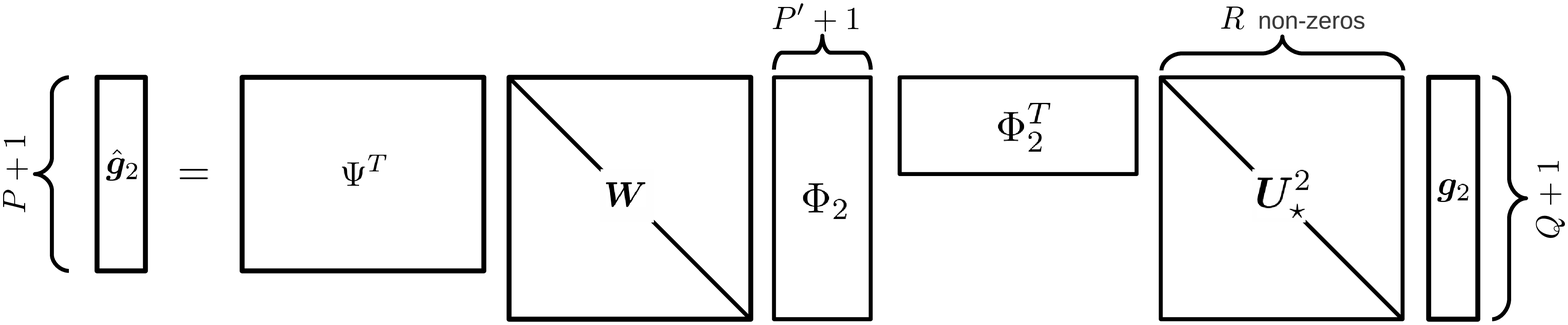}
        }\\
\subfloat{
\includegraphics[width=0.9\textwidth]{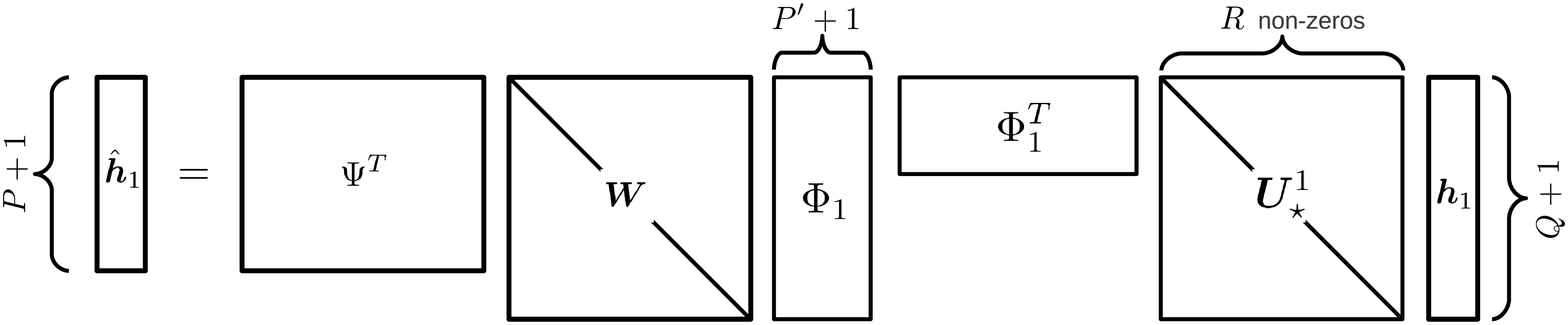}
        }\\
\subfloat{
\includegraphics[width=0.9\textwidth]{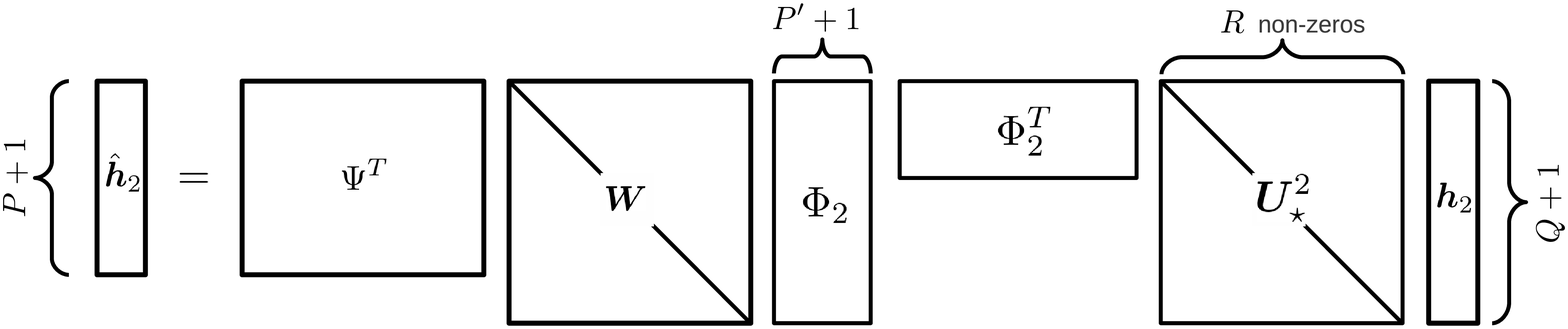}
        }
\caption{Matrix representation of pseudospectral coefficient approximation with
  reduced bases and modified quadrature for
  the vectors of coefficients $\hat{\bm{g}}_1$, $\hat{\bm{g}}_2$,
  $\hat{\bm{h}}_1$, and $\hat{\bm{h}}_2$.}
\label{fig:ps2}
\end{figure}

With reduced polynomial bases and modified weight sets in hand, we can
approximate the coefficients needed to construct the Galerkin
system~\eqref{eq:newtongalerkin} using $R_i<Q+1$ solves of each PDE system. More
precisely, we approximate the coefficients, $\hat{\bm{g}}_1$, $\hat{\bm{g}}_2$,
$\hat{\bh}_1$, and $\hat{\bh}_2$ using the computations written in matrix
notation in Figure \ref{fig:ps2}, where $\bm{U}_\star^1$ and $\bm{U}_\star^2$ are
diagonal matrices of the modified weight vectors $\bm{u}_\star^1$ and
$\bm{u}_\star^2$, respectively. 

\vspace{-6pt}


\section{Numerical examples} \label{sec:examples}

We now study the performance of the above approaches on two problems.  First we
apply the procedures for generating the reduced basis and associated modified
quadrature rule to a simple composite function where the accuracy of the
techniques can be easily explored.  We then apply them to a set of network
coupled PDEs representing an idealized model of a nuclear reactor,
demonstrating a significant reduction in computational cost.

\vspace{-2pt}

\subsection{Simple composite function}

As a test problem to demonstrate the accuracy of the basis reduction and
modified quadrature procedures, consider the following composite function:
\begin{equation}\label{eq:simple_func}
	\begin{split}
		y_1(x) &= x_1, \\
		y_2(x) &= \frac{1}{10 + \sum_{i=1}^{s} \frac{x_i}{i}}, \\
		h(y) &= \exp\left(y_1 + y_2\right),
	\end{split}
\end{equation}
for $x\in\R^s$ and $y=(y_1,y_2)\in\R^2$.  Here we consider $s=4$ and allow each
$x_i$ to be a uniform random variable over $[-1,1]$.  This form of composite
function, where some components of $y$ depend linearly on $x$, was chosen
because of its similarity to the structure of the network problem described
later.  Even so, the components of $y$ are clearly dependent.  For a given
polynomial order $N$, we build the pseudospectral approximation
$y(x)\approx\tilde{y}(x) = \sum_{k=0}^P \hat{y}_k\psi_k(x)$, where the polynomials
$\psi_k(x)$ are products of normalized Legendre polynomials of total order at
most $N$, using Gauss-Legendre tensor product quadrature with $Q+1 = (N+1)^s$
points.  We then construct pseudospectral approximations
$h(\tilde{y}(x))\approx\tilde{h}(x) = \sum_{k=0}^P \hat{h}_k\psi_k(x)$ by sampling
$h(\tilde{y}(x))$ on the same quadrature grid and using the reduced basis and
modified quadrature techniques described above.  In the latter case we build the
approximation~\eqref{eq:pc_y} with $N'=N$.  Define $\bm{\hat{h}}^{(N)}$ to be
the vector of pseudspectral coefficients in the first case and
$\bm{\tilde{h}}^{(N)}$ in the second.

In Table~\ref{tab:simple_d4} we vary $N$ from 1 to 10 and display the polynomial
basis size $P+1$, the quadrature size $Q+1$, the reduced basis size $P'+1$, the
number of non-zero weights $R$ in the modified quadrature, the largest error in
the pseudospectral coefficients, and largest error in the discrete orthogonality
of the reduced basis with the modified quadrature.  For each $N$, we used the
simplex solver in the Clp package~\cite{ClpURL} to solve the linear
program~\eqref{eq:lin_prog_final}.  The error in the pseudospectral coefficients
is computed by comparing them to the $N=10$ solution using the original
(unreduced) basis and and quadrature.  We use a tolerance of
$10^{-12}$ in the QR factorization~\eqref{eq:red_quad_qr}.  One
can see a dramatic reduction in the number of samples of $h$ needed
(approximately an order of magnitude), with essentially no difference in the
error in the coefficients.  In Table~\ref{tab:simple_d4_loose} we repeat the
experiments with a tolerance of only $10^{-6}$ in the QR factorization.  Here
one can see this reduces the number of non-zero quadrature weights even more,
but dramatically increases the discrete orthogonality error.  This additional
error pollutes the pseudospectral coefficients once the coefficient error
becomes roughly the same order of magnitude as the discrete orthogonality error.

The matrix $\mA$ used to generate the modified quadrature
constraint~\eqref{eq:red_constraint} should in exact arithmetic have rank $2N'+2
\choose 2$ for this problem, but finite precision computations destroys this.
To demonstrate this, we repeat the experiments in
Table~\ref{tab:simple_d4_restrict} taking the first $2N'+2 \choose 2$ columns of
$\mQ$ in~\eqref{eq:red_quad_qr} instead of using a tolerance on $\mS$.  Note that
this is similar to an approach investigated previously~\cite{Arnst:2012cn} where
the constraint was based on generating the reduced basis polynomials up to order
$2N'$.  We see performance similar to using a looser tolerance in the QR
factorization, namely fewer non-zero quadrature weights at the expense of
increased discrete orthogonality error that at higher order reduces the accuracy
of the coefficients.  Thus basing the modified quadrature linear program on
directly maintaining discrete orthogonality~\eqref{eq:lin_prog}, and employing a
relatively tight tolerance to extract a linearly independent set of constraints,
appears to be more robust.
\begin{table}[htbp]
  \centering
  \caption{Performance of the reduced basis and modified quadrature approaches with a tolerance of $10^{-12}$ on the QR factorization of the linear program constraint.}
    \begin{tabular}{cccccccc}
    \toprule
    $N$ & $P+1$ & $Q+1$ & $P'+1$ & $R$ & $\|\bm{\hat{h}}^{(10)}-\bm{\hat{h}}^{(N)}\|_\infty$ & $\|\bm{\hat{h}}^{(10)}-\bm{\tilde{h}}^{(N)}\|_\infty$ & $\|\mbox{vec}(\bm{I} - \bm{Z}^T \bm{U} \bm{Z})\|_\infty $ \\
    \midrule
    1     & 5     & 16    & 3     & 5     & 2.93E-02 & 2.93E-02 & 2.22E-16 \\
    2     & 15    & 81    & 6     & 12    & 3.58E-03 & 3.58E-03 & 2.10E-14 \\
    3     & 35    & 256   & 10    & 22    & 3.55E-04 & 3.55E-04 & 1.46E-12 \\
    4     & 70    & 625   & 15    & 47    & 2.94E-05 & 2.94E-05 & 1.37E-12 \\
    5     & 126   & 1296  & 21    & 101   & 2.09E-06 & 2.09E-06 & 1.83E-12 \\
    6     & 210   & 2401  & 28    & 188   & 1.30E-07 & 1.30E-07 & 2.55E-12 \\
    7     & 330   & 4096  & 36    & 346   & 7.18E-09 & 7.18E-09 & 3.81E-12 \\
    8     & 495   & 6561  & 45    & 587   & 3.58E-10 & 3.57E-10 & 6.10E-12 \\
    9     & 715   & 10000 & 55    & 941   & 1.62E-11 & 1.63E-11 & 2.62E-12 \\
    10    & 1001  & 14641 & 66    & 1425  & 0.00E+00 & 1.63E-12 & 2.70E-12 \\
    \bottomrule
    \end{tabular}%
  \label{tab:simple_d4}%
\end{table}%
\begin{table}[htbp]
  \centering
   \caption{Performance of the reduced basis and modified quadrature approaches with a tolerance of $10^{-6}$ on the QR factorization of the linear program constraint.}
    \begin{tabular}{cccccccc}
    \toprule
    $N$ & $P+1$ & $Q+1$ & $P'+1$ & $R$ & $\|\bm{\hat{h}}^{(10)}-\bm{\hat{h}}^{(N)}\|_\infty$ & $\|\bm{\hat{h}}^{(10)}-\bm{\tilde{h}}^{(N)}\|_\infty$ & $\|\mbox{vec}(\bm{I} - \bm{Z}^T \bm{U} \bm{Z})\|_\infty $ \\
    \midrule
    1     & 5     & 16    & 3     & 5     & 2.93E-02 & 2.93E-02 & 2.22E-16 \\
    2     & 15    & 81    & 6     & 12    & 3.58E-03 & 3.58E-03 & 2.10E-14 \\
    3     & 35    & 256   & 10    & 22    & 3.55E-04 & 3.55E-04 & 1.46E-12 \\
    4     & 70    & 625   & 15    & 35    & 2.94E-05 & 2.94E-05 & 1.90E-11 \\
    5     & 126   & 1296  & 21    & 50    & 2.09E-06 & 1.83E-06 & 2.28E-06 \\
    6     & 210   & 2401  & 28    & 70    & 1.30E-07 & 1.22E-07 & 4.02E-08 \\
    7     & 330   & 4096  & 36    & 92    & 7.18E-09 & 4.24E-07 & 1.11E-06 \\
    8     & 495   & 6561  & 45    & 158   & 3.58E-10 & 4.23E-06 & 1.12E-03 \\
    9     & 715   & 10000 & 55    & 252   & 1.62E-11 & 3.86E-06 & 4.87E-06 \\
    10    & 1001  & 14641 & 66    & 475   & 0.00E+00 & 1.30E-05 & 3.85E-02 \\
    \bottomrule
    \end{tabular}%
     \label{tab:simple_d4_loose}%
\end{table}%

\begin{table}[htbp]
  \centering
  \caption{Performance of the reduced basis and modified quadrature approaches using the first $2N'+2 \choose 2$ columns in the QR factorization of the linear program constraint.}
    \begin{tabular}{cccccccc}
    \toprule
    $N$ & $P+1$ & $Q+1$ & $P'+1$ & $R$ & $\|\bm{\hat{h}}^{(10)}-\bm{\hat{h}}^{(N)}\|_\infty$ & $\|\bm{\hat{h}}^{(10)}-\bm{\tilde{h}}^{(N)}\|_\infty$ & $\|\mbox{vec}(\bm{I} - \bm{Z}^T \bm{U} \bm{Z})\|_\infty $ \\
    \midrule
    1     & 5     & 16    & 3     & 5     & 2.93E-02 & 2.93E-02 & 2.22E-16 \\
    2     & 15    & 81    & 6     & 12    & 3.58E-03 & 3.58E-03 & 2.10E-14 \\
    3     & 35    & 256   & 10    & 22    & 3.55E-04 & 3.55E-04 & 1.46E-12 \\
    4     & 70    & 625   & 15    & 45    & 2.94E-05 & 2.94E-05 & 1.97E-11 \\
    5     & 126   & 1296  & 21    & 66    & 2.09E-06 & 2.09E-06 & 2.19E-10 \\
    6     & 210   & 2401  & 28    & 91    & 1.30E-07 & 1.40E-07 & 2.76E-08 \\
    7     & 330   & 4096  & 36    & 120   & 7.18E-09 & 9.16E-08 & 2.46E-07 \\
    8     & 495   & 6561  & 45    & 153   & 3.58E-10 & 5.34E-06 & 1.05E-05 \\
    9     & 715   & 10000 & 55    & 190   & 1.62E-11 & 1.46E-04 & 4.24E-04 \\
    10    & 1001  & 14641 & 66    & 231   & 0.00E+00 & 3.63E-03 & 6.06E-03 \\
    \bottomrule
    \end{tabular}%
  \label{tab:simple_d4_restrict}%
\end{table}%

\subsection{Nuclear reactor simulation}
\begin{figure}
	\centering
	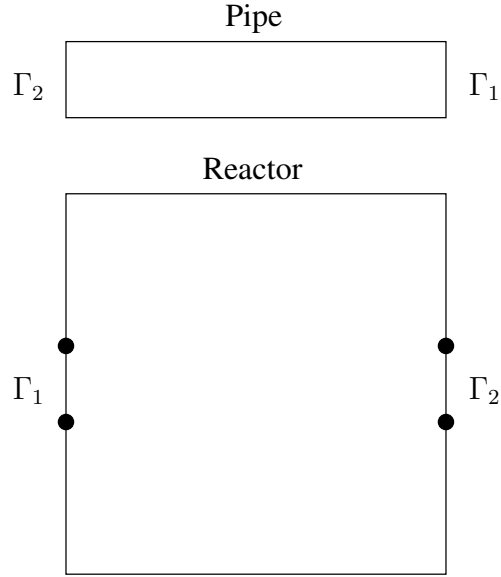
	\caption{Pipe and reactor geometries.  Fluid flows from the pipe into the reactor across $\Gamma_1$ and out of the reactor into the pipe across $\Gamma2$.}
\label{fig:reactor}
\end{figure}

Finally, we consider a coupled PDE problem, inspired from uncertainty
quantification of nuclear reactor simulations, and demonstrate improvement in
total computational cost using the reduced basis and quadrature techniques
discussed above.  Typically these models consist of a high fidelity model of the
nuclear reactor core coupled to a low dimensional network model of the rest of
the plant representing the various pipes, heat exchanges, turbines, and so on.
Uncertainties can arise in any component of this model, and if the system is
strongly coupled, these uncertainties must be tracked throughout the whole
system.  Since there can be many components in the network model giving rise to
a very large number of independent sources of uncertainty, and since the reactor
core simulations themselves are typically quite computationally expensive, it is
often infeasible to propagate all of the uncertainties that arise in the coupled
plant-core model and accurately resolve their effects.  Thus it becomes
necessary to propagate uncertainties in each component (such as the core)
separately, making assumptions on how these uncertainties interact with the rest
of the system.  However by using the basis and quadrature reduction techniques
detailed above, we can make the computational cost of expensive components in
the model (again, such as the core) effectively independent of the number of
independent sources of uncertainty arising from other components in the model.
To demonstrate this, we consider a model consisting of two coupled
thermal-hydraulics components: a $1\times 0.1$ in/out-flow pipe and a $1\times
1$ reactor vessel.  Fluid flows in to the reactor from the pipe, is heated, and
flows out of the reactor back into the pipe where it is cooled.  The
computational domain is shown in Figure~\ref{fig:reactor}.
We include a temperature source in the reactor (to represent heating of the
fluid from fission) and cool the upper and lower walls of the pipe by holding
them at a fixed temperature (which conceptualizes the rest of the
thermal-hydraulics of the plant).  In each component we solve the coupled
steady-state Navier Stokes and energy equations,
\begin{equation}\label{eq:navier_stokes}
	\begin{split}
		-\nu \Delta \bu + \bu \cdot \nabla \bu + \nabla p =& \beta(T-T_{\text{ref}}) \bm{g}, \\
		-\kappa \Delta T + \bu \cdot \nabla T + T_s =& 0, 
	\end{split}
\end{equation}
where $\nu$ is the kinematic viscosity, $\rho$ is the density, $\beta$ is the
coefficient of thermal expansion, $\bm{g}$ is the gravity vector, $\kappa$ is
the thermal diffusivity and $T_s$ is the heat source.  We let $\Gamma_{1}$ and
$\Gamma_{2}$ denote the respective interfaces between the inflow of the pipe to
the reactor, and the outflow of the reactor to the pipe.  We use $\bu_{P}$ and
$\bu_{R}$ to denote the fluid velocities in the pipe and reactor, respectively.
We use similar subscript notation for the pressure and temperature fields.
Apart from the interfaces $\Gamma_1$ and $\Gamma_2$, the fluid velocity and
temperature are fixed at all boundaries to zero.

The fluid variables are coupled through the following interface conditions,
\begin{equation}\label{eq:fluid_interface1}
 \begin{cases}
  \bu_{P} = \bu_{R}, & x\in \Gamma_{1}, \\
  (\nu \nabla \bu_{P} - p_{P} \mathbb{I})\bm{n}_1 = (\nu \nabla \bu_{R} - p_R \mathbb{I})\bm{n}_1, & x\in \Gamma_{1},
 \end{cases}
\end{equation}
and
\begin{equation}\label{eq:fluid_interface2}
 \begin{cases}
  \bu_{R} = \bu_{P}, & x\in \Gamma_{2}, \\
  (\nu \nabla \bu_{R} - p_{R} \mathbb{I})\bm{n}_2 = (\nu \nabla \bu_{P} - p_{P} \mathbb{I})\bm{n}_2, & x\in \Gamma_{2},
 \end{cases}
\end{equation}
where $\bm{n}_1$ and $\bm{n}_2$ are the unit outward normal to the reactor on
$\Gamma_1$ and $\Gamma_2$ respectively.
The temperature fields are coupled through similar interface conditions,
\begin{equation}\label{eq:temp_interface1}
 \begin{cases}
  T_{P} = T_{R}, & x\in \Gamma_{1}, \\
  \kappa \nabla T_{P} \cdot \bm{n}_1 = \kappa \nabla T_{R} \cdot \bm{n}_1, & x\in \Gamma_{1},
 \end{cases}
\end{equation}
and
\begin{equation}\label{eq:temp_interface2}
 \begin{cases}
  T_{R} = T_{P}, & x\in \Gamma_{2}, \\
  \kappa \nabla T_{R} \cdot \bm{n}_2 = \kappa \nabla T_{P} \cdot \bm{n}_2, & x\in \Gamma_{2}.
 \end{cases}
\end{equation}
We recall the non-dimensional Reynolds number,
\[\text{Re} = \frac{UL}{\nu},\]
where $U$ and $L$ are reference velocities and length respectively and the
nondimensional Rayleigh number,
\[\text{Ra} = \frac{|\bm{g}|\beta (T_1-T_2) L^3}{\nu \kappa},\]
where $(T_1-T_2)>0$ is a reference temperature difference.  A weak solution to
the coupled system can be shown to exist under certain assumptions on the
Reynolds and Rayleigh numbers \cite{CKN99,CI99}.

The components outside of the reactor core are typically modeled at much lower
fidelity than the core itself, and are coupled to the core through low
dimensional interfaces which define the network model.  Thus we define a network
problem by simplifying the interface conditions so that only the spatial
averages are coupled at the interfaces.  We use $\bbu$, $\bbp$ and $\bbT$ (with
the appropriate subscripts) to denote the spatial average of these fields over
one of the interfaces.  We do not explicitly specify which interface since this
is clear from the context.  In the pipe and reactor, the chosen boundary
conditions imply conservation of mass guaranteeing the average flow coming into
the pipe/reactor is equal to the average flow leaving.  Thus, the {\em fluid}
interface conditions on $\Gamma_1$~\eqref{eq:fluid_interface1} and
$\Gamma_2$~\eqref{eq:fluid_interface2} are automatically satisfied on average.
Therefore, the interface coupling conditions that must be explicitly enforced
are
\begin{equation}\label{eq:avgtemp_interface1}
 \begin{cases}
  \bbT_{P} = \bbT_{R}, & x\in \Gamma_{1}, \\
  \kappa \nabla \bbT_{P} \cdot \bm{n}_1 = \kappa \nabla \bbT_{R} \cdot \bm{n}_1, & x\in \Gamma_{1},
 \end{cases}
\end{equation}
and
\begin{equation}\label{eq:avgtemp_interface2}
 \begin{cases}
  \bbT_{R} = \bbT_{P}, & x\in \Gamma_{2}, \\
  \kappa \nabla \bbT_{R} \cdot \bm{n}_2 = \kappa \nabla \bbT_{P} \cdot \bm{n}_2, & x\in \Gamma_{2}.
 \end{cases}
\end{equation}

For each of the reactor and pipe components, the Navier-Stokes
equations~\eqref{eq:navier_stokes} are discretized in space using bi-linear
finite elements with pressure stabilization (PSPG) and upwinding
(SUPG)~\cite{SHADID:2006km}.  The reactor and pipe geometries are represented by
quadrilateral meshes with side length $1/40$ (1600 and 160 mesh cells
respectively).  We treat the thermal diffusivity of the pipe as uncertain, by
representing it as a correlated random field with exponential covariance (with
mean 0.1, standard deviation 0.05, and correlation lengths of 0.1 and 0.01 in
the flow and transverse directions respectively) and approximated via a
truncated Karhunen-Lo\`{e}ve expansion in $s$ terms~\cite{Ghanem_Spanos_91}.  We
assume the resulting random variables are uniformly distributed on $[-1,1]$ and
are independent.  The spatially discretized equations in each component along
with the discretized interface
conditions~\eqref{eq:avgtemp_interface1},~\eqref{eq:avgtemp_interface2} can then
be written as the abstract coupled system
\begin{equation}\label{eq:pipe_reactor_network}
  \begin{split}
    f_1(u_1,v_1,v_2,\xi) &= 0, \\
    f_2(u_2,v_1,v_2) &= 0, \\
    g_1(u_1) - g_{2}(u_2) &= 0, \\
    g_{3}(u_2) - g_4(u_1) &= 0,
 \end{split}
\end{equation}
where $f_1$ and $f_2$ are the finite element residual equations for the pipe and
reactor, $u_1$ and $u_2$ are the nodal vectors of fluid velocity and temperature
unknowns in the pipe and reactor, $v_1$ and $v_2$ represent the Neumann data
($\kappa \nabla \bbT \cdot \bm{n}$) on $\Gamma_1$ and $\Gamma_2$, $\xi =
(\xi_1,\dots,\xi_s)$ are the Karhunen-Lo\`{e}ve random variables, and
\begin{equation}\label{eq:pipe_reactor_network_funcs}
\begin{split}
g_1(u_1) &= \bbT_{P}|_{\Gamma_1}, \\
g_{2}(u_2) &= \bbT_{R}|_{\Gamma_1}, \\ 
g_{3}(u_2) &= \bbT_{R}|_{\Gamma_2}, \\
g_4(u_1) &= \bbT_{P}|_{\Gamma_2}. \\
\end{split}
\end{equation}
Through elimination of $u_1$ and $u_2$ for realizations of $\xi$ these become
Neumann-to-Dirichlet maps that are well-posed due to our choice of boundary
conditions away from the interfaces.  For a deterministic problem, this is a
straightforward modification of a nonoverlapping domain decomposition method
\cite{QuarterVal,DDBook}.

\begin{figure}{t}
\centering
\subfloat[Reactor temperature mean]{
\includegraphics[width=0.45\textwidth]{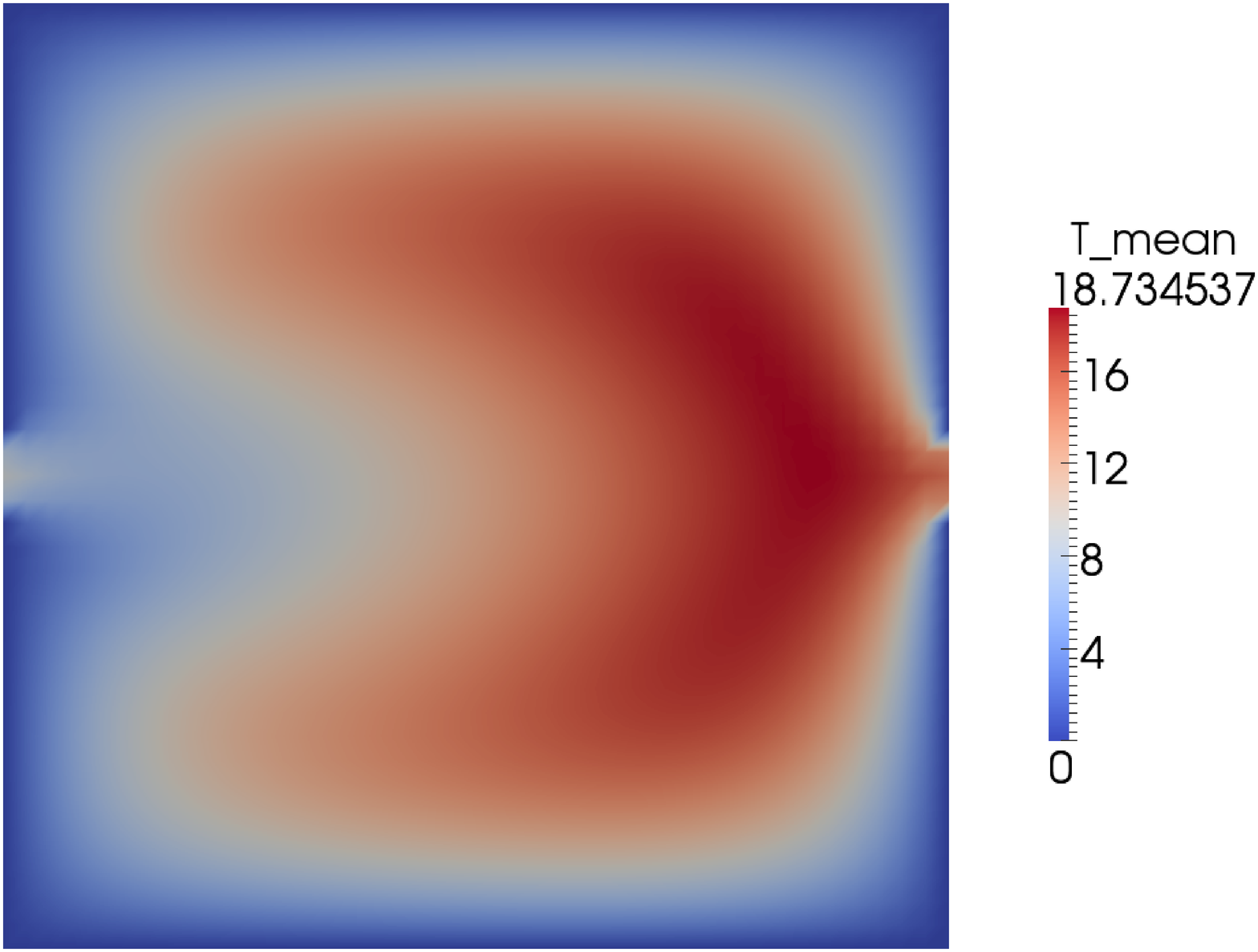}
\label{fig:reactor_mean}
}\quad
\subfloat[Reactor temperature std. dev.]{
\includegraphics[width=0.45\textwidth]{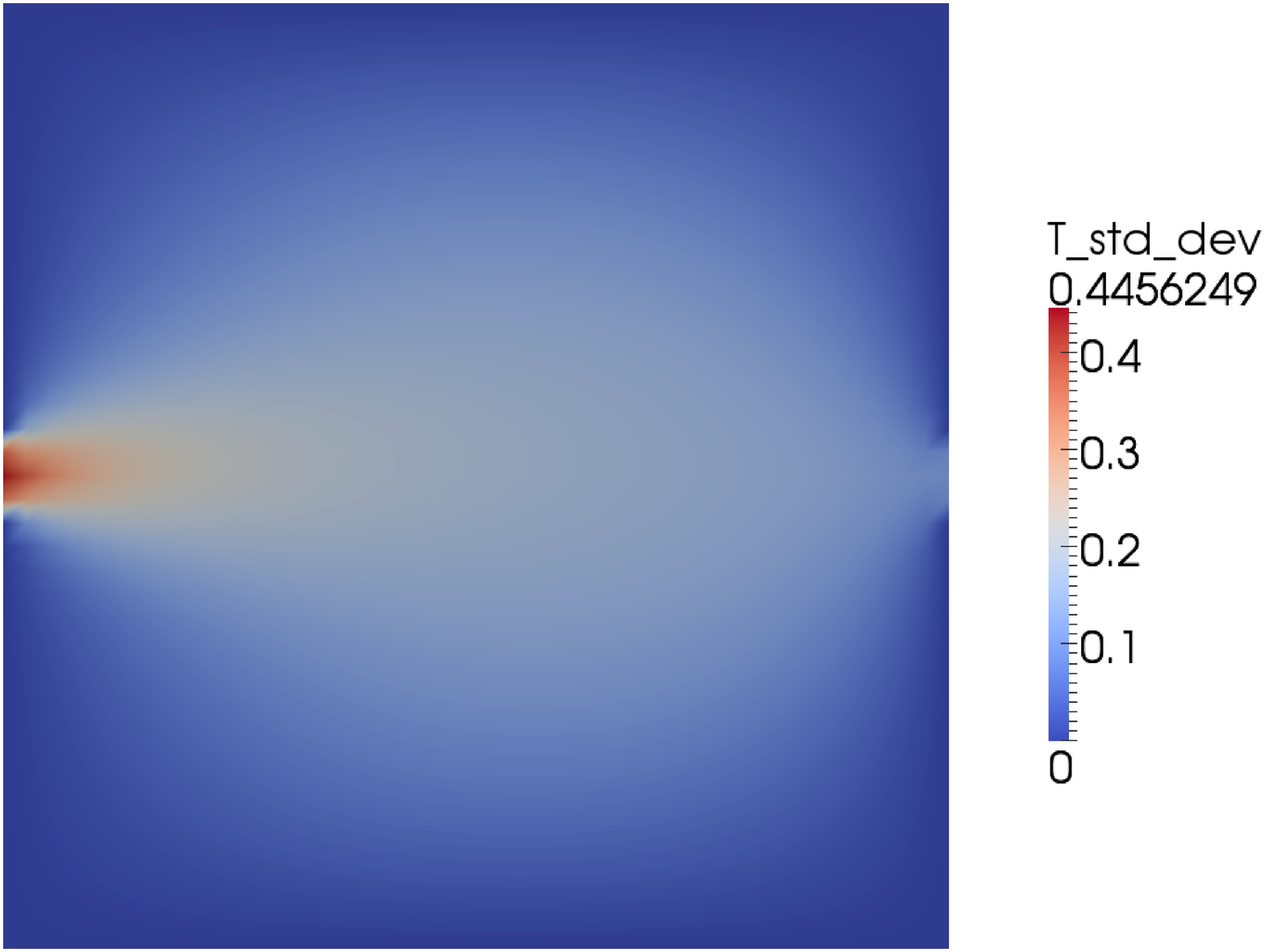}
\label{fig:reactor_sd}
}\\
\subfloat[Pipe temperature mean]{
\includegraphics[width=0.45\textwidth]{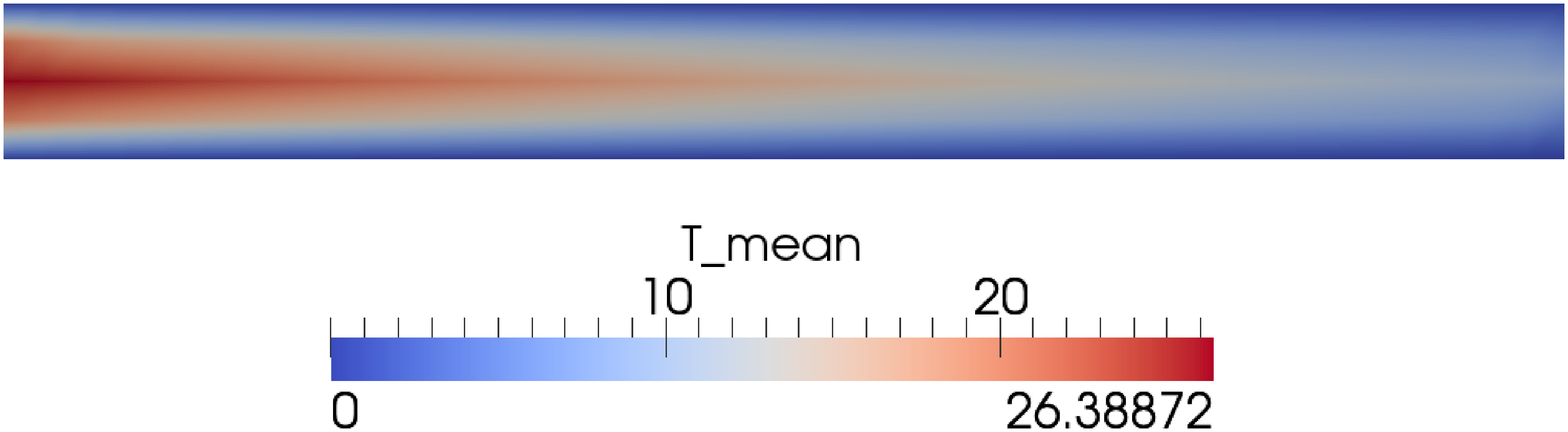}
\label{fig:pipe_mean}
}\quad
\subfloat[Pipe temperature std. dev.]{
\includegraphics[width=0.45\textwidth]{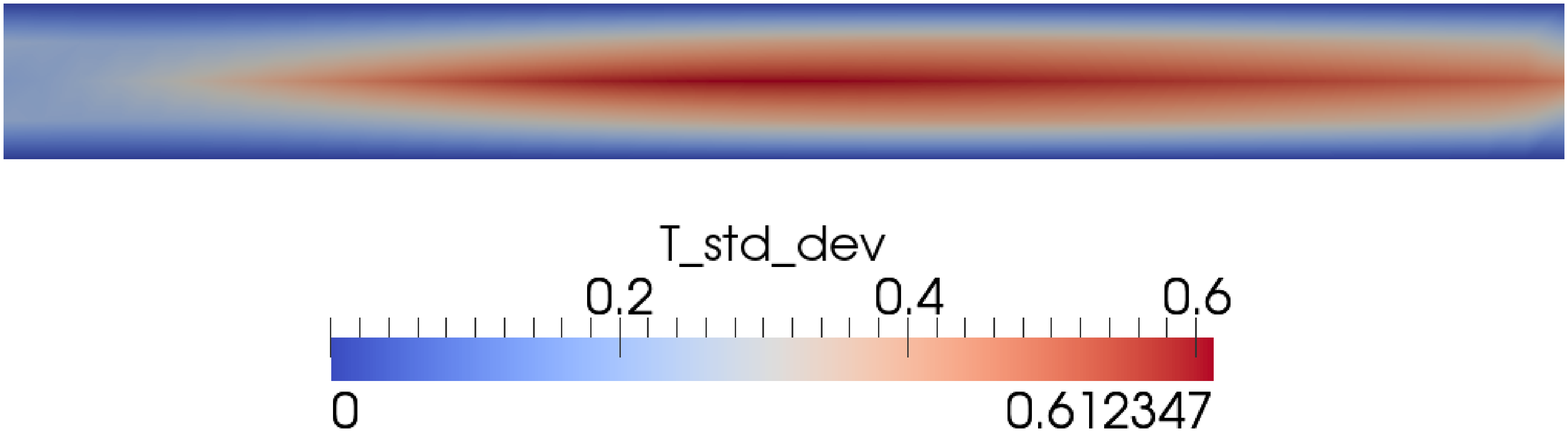}
\label{fig:pipe_sd}
}
\caption{Reactor and pipe temperature mean and standard deviation for $s=3$ (not drawn to scale).}
\label{fig:temp}
\end{figure}
This form of network coupled system is a slight generalization of the
system~\eqref{eq:stoch_general_multiphysics}, and the techniques presented above
can be easily extended to this problem.  Given the number of random variables
$s$, we apply the nonlinear elimination-based stochastic Galerkin procedure as
outlined in Section~\ref{sec:multi_phys_galerkin} using a total polynomial order
$N=3$.  We compute the pseudospectral approximation (see~\Eqref{eq:ps}) of the
pipe and reactor coupling functions~\eqref{eq:pipe_reactor_network_funcs} using
a tensor-product Gauss-Legendre grid with $N+1$ points for each random variable.
Even though uncertainty only appears in the pipe component, its effects must be
propagated through the reactor due to the coupling variables $v_1$ and $v_2$,
thus resulting in numerous solves of the reactor problem $f_2=0$.  All
calculations were implemented within the Albany simulation environment
leveraging numerous Trilinos~\cite{TrilinosTOMS} packages for assembling and
solving the pipe and reactor finite element equations.  The finite element
solution for each realization of the random variables was computed via a
standard Newton method employing a GMRES iterative linear solver preconditioned
by an incomplete LU factorization preconditioner.  The resulting stochastic
Galerkin network equations were assembled and solved using the Stokhos intrusive
stochastic Galerkin package~\cite{StokhosURL}.  Since the network system after
nonlinear elimination is just a $2\times 2$ system, the cost of solving this
system is insignificant compared to the total cost of solving the reactor
component at each realization of $\xi$.  Figure~\ref{fig:temp} shows the
resulting mean and standard deviation of the fluid temperature for the uncertain
diffusivity field with coefficient of variation of $0.5$ and $s=3$ random
variables.  The reactor has a temperature source $T_s=500$, Reynolds number
$\text{Re} = 100$ and mean Rayleigh number $\text{Ra}\approx 0.87$.

We then apply the basis reduction and quadrature methods described in
Section~\ref{sec:transformation} to the reactor component only.
Table~\ref{tab:reactor} compares the resulting total run time for the complete
stochastic Galerkin, nonlinear elimination network calculation between both
approaches for an increasing number of random variables in the pipe.  One can
see that for the reduced basis and quadrature approach, the size of the
polynomial basis, modified quadrature, and resulting solution time for the
reactor component is roughly constant, even though the total number of random
variables appearing in the coupled system is increasing.
\begin{table}[tbp]
  \centering
  \caption{Performance comparison of stochastic Galerkin applied to the pipe-reactor network system~\eqref{eq:pipe_reactor_network} using the traditional and reduced basis/quadrature approaches.}
    \begin{tabular}{ccccccccccc}
    \toprule
    & & & & & \multicolumn{3}{c}{Time (sec)} &   \multicolumn{3}{c}{Reduced Time(sec)} \\ \cmidrule(lr){6-8} \cmidrule(lr){9-11}
    $s$ & $P+1$ & $Q+1$ & $P'+1$ & $R$ & Pipe & Reactor & Total & Pipe & Reactor & Total \\
    \midrule
    2 &            10 &       16 &        10 &        16 &        4 &        62 &          67 &     4 &         53 &       58 \\
    3 &            20 &       64 &        10 &        40 &      17 &      246 &       263 &    17 &      120 &     137 \\
    4 &            35 &     256 &        10 &        41 &      82 &   1052 &     1134 &    73 &      129 &     202 \\
    5 &            56 &   1024 &        10 &        35 &    353 &   4051 &     4405 & 341 &      116 &     458 \\
    \bottomrule
    \end{tabular}%
  \label{tab:reactor}%
\end{table}%

\vspace{-6pt}

\section{Summary \& Conclusions}\label{sec:conclusions}

We have presented a method for constructing a polynomial surrogate response
surface for the outputs of a network coupled multiphysics system that exploits
its structure to increase efficiency. We reduce the full system with a nonlinear
elimination method, which results in a smaller system to solve for the coupling
terms that depends on the uncertain inputs represented as parameters. We then
apply a stochastic Galerkin procedure with a Newton iteration to compute the
coefficients of a surrogate response surface that approximates the coupling
terms as a polynomial of the input parameters. The residual and Jacobian matrix
in the Newton update can be viewed as composite functions: these terms depend on
the coupling terms which depend on the uncertain inputs. We take advantage of
this composite structure to build a reduced polynomial basis that depend on the
coupling terms as intermediate variables. We use this reduced basis to find a
modified quadrature rule with relatively few nonzero weights. Each weight equal
to zero corresponds to a PDE solve that can be ignored when solving the Newton
system. This results in substantial computational savings, which we demonstrated
on a simple model of a nuclear reactor.

The method is appropriate when the number of input parameters to the
full system is small enough to work with a tensor product quadrature
grid. Even though we only evaluate the relatively cheap coupling terms
at the full grid to construct the modified quadrature rule, we still
have an exponential increase in the number of points in the grid as
the dimension of the input space increases; a cheap function evaluated
$10^{10}$ times can still be expensive. This can be alleviated to some
extent by standard methods for sensitivity analysis and anisotropic
approximation. In principle, sparse grids could also be used in the
method by relaxing the non-negativity constraint in the linear
program~\eqref{eq:lin_prog_final}, but one must take great care when
using sparse grids in conjunction with integration for pseudospectral
approximation to maintain discrete orthogonality~\cite{Spam:2012}, as
well as insuring positive definiteness of the inner products implicit
in the QR factorizations~\eqref{eq:basis_qr}
and~\eqref{eq:red_quad_qr}.  This last difficulty could be alleviated
for computing the reduced basis by projecting the monomial matrix
$\mY$~\eqref{eq:monomials} onto the original polynomial basis
$\{\psi_i\}$, and thus~\eqref{eq:basis_qr} becomes a standard
(unweighted) QR factorization.  However one must still choose a sparse
grid when defining the modified quadrature rule that has degree of
exactness at least $2N'$ for~\eqref{eq:red_quad_qr}.  Alternatively,
in the spirit of~\cite{Spam:2012}, one could incorporate this into a
Smolyak sparse grid approach by applying the reduced basis and
quadrature techniques to each tensor grid appearing in the Smolyak
expansion.

\bibliographystyle{wileyj}
\bibliography{paper}
\end{document}